\documentclass[a4paper,times,3p]{elsarticle}

\usepackage{lineno,hyperref}
\usepackage{setspace}

\usepackage{amsmath}
\usepackage[table,xcdraw]{xcolor}
\usepackage{tabularx}
\usepackage{multirow}
\usepackage{subcaption} 
\usepackage{csvsimple}
\usepackage{caption}
\usepackage{array,multirow}
\usepackage{todonotes}
\usepackage{tikz}
\usepackage{pgfplots}
\usepackage{floatrow}
\usepackage{tabto}
\usepackage[per-mode=symbol]{siunitx}
\newfloatcommand{capbtabbox}{table}[][\FBwidth]

\usepackage{blindtext}

\usepackage{makecell}
\usepackage{xstring}
\usepackage{comment}
\usepackage{cleveref}

\pgfplotsset{compat = newest}
\usepgfplotslibrary{units}

\usepackage[]{algorithm2e}

\newcommand{\qty}[1]{\SI{#1}}

\newcommand{\vekt}[1]{\mathbf{#1}}
\newcommand{\fvel}[0]{\vekt{v}^f}
\newcommand{\CauchyStress}[0]{\vekt{T}}
\newcommand{\nsd}[0]{{n_{sd}}}
\newcommand{\foralltime}[0]{\forall t \ge 0}

\newcommand{\Rk}[1]{\vect{r}^{#1}}
\newcommand{\Jac}[0] { \widehat{\vekt{J}^{-1}}}

\newcommand{\Equ}[1]{Eq.~(\ref{equ:#1})}
\newcommand{\Eqs}[1]{Eqs.~(\ref{eqs:#1})}
\newcommand{\Fig}[1]{Figure~\ref{fig:#1}}
\newcommand{\Tab}[1]{Table~\ref{tab:#1}}

\newcommand{\Sec}[1]{Section~\ref{sec:#1}}
\newcommand{\App}[1]{\ref{app:#1}}

\newcommand{\Rem}[1]{Remark~\ref{Remark:#1}}
\newcommand{\Fnt}[1]{footnote~\ref{fnt:#1}}

\newcommand{\uv}[0]{\vect{u}}
\newcommand{\bv}[0]{\vect{b}}
\newcommand{\Am}[0]{\mat{A}}
\newcommand{\Mm}[0]{\mat{M}}
\newcommand{\resSolid}[0]{\vekt{r}_s}
\newcommand{\resFlow}[0]{\vekt{r}_f}
\newcommand{\resProblem}[0]{\vekt{r}_p}

\newcommand{\abs}[1]{\left\lvert#1\right\rvert}
\newcommand{\norm}[1]{\left\lVert#1\right\rVert}
\newcommand{\normE}[1]{\norm{#1}_2}

\newcommand{\flowDomain}[0]{\Omega^f}
\newcommand{\solidDomain}[0]{\Omega^s}

\newcommand{\minval}{0.0} 
\newcommand{\maxval}{1.0}
\definecolor{high}{HTML}{ff0000}
\definecolor{mid}{HTML}{ffff00}
\definecolor{low}{HTML}{22b20c}
\newcommand{\opacity}{40}
\newcommand{\ApplyGradient}[1]{
    \pgfmathsetmacro{\midval}{(\minval+\maxval)/2}
    \IfStrEq{#1}{-}{
        \cellcolor{high!\opacity} #1
    }{
        \IfDecimal{#1}{
            \ifdim #1 pt > \midval pt
                \pgfmathparse{max(min(100*(#1-\midval)/(\maxval-\midval),100),0)}
                \xdef\tempa{\pgfmathresult}
                \hspace{-0.33em}\cellcolor{high!\tempa!mid!\opacity} #1
            \else
                \pgfmathparse{max(min(100*(\midval-#1)/(\midval-\minval),100),0)}
                \xdef\tempa{\pgfmathresult}
                \hspace{-0.33em}\cellcolor{low!\tempa!mid!\opacity} #1
            \fi
        }{
            #1
        }
    }
}
\newcommand{\eqTime}[1]{\textbf{\ApplyGradient{#1}}}
\newcommand{\cIt}[1]{\textit{\textcolor{black}{#1}}}
\newcommand{\fIt}[1]{\textcolor{blue}{#1}}
\newcommand{\sIt}[1]{\textcolor{orange}{#1}}

\newcommand{\subproblemIter}[0]{subproblem iteration}

\newcommand{\mat} [1]{\mathbf{#1}}
\newcommand{\vect}[1]{\mathbf{#1}}
\newcommand{\E}[1]{\num{e#1}}

\newcommand{\inR}[1]{\in \mathbb{R}^{#1}}
\newcommand{\ofxt}[0]{(\vect{x},t)}
\newcommand{\nInterfaceDofs}[0]{n_\Gamma}

\newcommand{\costFac}[1]{C_{\text{#1}}}
\newcommand{\costCouple}[0]{\costFac{}^c}
\newcommand{\costSolver}[0]{\costFac{}}
\newcommand{\costIter}[0]{\costFac{iter}}
\newcommand{\costFix}[0]{\costFac{fix}}
\newcommand{\costSimulation}[0]{\costFac{simulation}}
\newcommand{\costCoupleSum}[0]{\overline{\costFac{}}^c}

\newcommand{\timeSimulation}[0]{t_{\text{simulation}}}

\newcommand{\totalCoupleIter}[0]{N^{c}}
\newcommand{\totalProblemIter}[0]{N}
\newcommand{\iterPerCall}[0]{n_{max}}

\newcommand{\kbar}[0]{\bar{k}}

\newcommand{\timing}[0]{t}

\newcommand{\numCalc}[0]{m}

\newcommand{\sumProblems}[0]{\sum_{p=f,s}}
\newcommand{\sumCplIter}[0]{\sum_{\kbar=1}^{\totalCoupleIter}}

\newcounter{para}
\newcommand \mypara{\refstepcounter{para} \par \noindent\textit{Remark \thepara:\space}}

\newcommand{\marktext}[1]{\textit{#1}}

\makeatletter
\def\blfootnote{\xdef\@thefnmark{}\@footnotetext}
\makeatother

\newcommand{\smallIfElsevier}[0]{\small}

\modulolinenumbers[5]
\graphicspath{ {Figures/} {../}}

\journal{International Journal for Numerical Methods in Engineering}

\biboptions{sort&compress}

\begin{document}

\begin{frontmatter}

\title{On the number of subproblem iterations per coupling step in partitioned fluid-structure interaction simulations}

\author[cats]{Thomas Spenke\fnref{cor}\corref{contri}} 
\ead{spenke@cats.rwth-aachen.de}

\author[ugent]{Nicolas Delaiss\'e\corref{contri}}
\ead{nicolas.delaisse@ugent.be}

\author[ugent,fm]{Joris Degroote}
\ead{joris.degroote@ugent.be}

\author[cats]{Norbert Hosters}
\ead{hosters@cats.rwth-aachen.de}

\fntext[cor]{Corresponding author}
\cortext[contri]{Contributed equally and should be considered joint first author}

\address[cats]{Chair for Computational Analysis of Technical Systems (CATS),\\Center for Simulation and Data Science (JARA-CSD),\\RWTH Aachen University,\\Schinkelstra\ss e 2, Aachen, Germany} 

\address[ugent]{Department of Electromechanical, Systems and Metal Engineering,\\Ghent University,\\Sint-Pietersnieuwstraat 41, Ghent, Belgium}

\address[fm]{Flanders Make @ UGent -- Core Lab MIRO}

\begin{abstract}
	In literature, the cost of a partitioned fluid-structure interaction scheme is typically assessed by the number of coupling iterations required per time step, while ignoring the internal iterations within the nonlinear subproblems.
	In this work, we demonstrate that these internal iterations have a significant influence on the computational cost of the coupled simulation.
    Particular attention is paid to how limiting the number of iterations within each solver call
    can shorten the overall run time, as it avoids polishing the subproblem solution using
    unconverged coupling data.
%
    Based on systematic parameter studies, we investigate the optimal number of \subproblemIter s per coupling step.
    %

	Lastly, this work proposes a new convergence criterion for coupled systems that is based on the residuals of the subproblems and therefore does not require any additional convergence tolerance for the coupling loop.
\end{abstract}

	


\end{frontmatter}


\section{Introduction}
\label{sec:introduction}

For the solution of fluid-structure interaction (FSI) problems, partitioned approaches are widespread, because they allow to reuse mature and reliable solvers, tailored to each of the two subproblems \cite{felippa2001partitioned}.
The flow and structure solvers are treated as black boxes and all data exchange is limited to the shared interface.
The main drawback of the partitioned technique is its need for an iterative coupling loop that repeatedly solves the subproblems within each time step to assure satisfaction of the equilibrium conditions on the interface.
Moreover, these coupling iterations are prone to stability issues, due to the added-mass effect \cite{Causin2005, Forster2007, VanBrummelen2009}.

During the last two decades, these instabilities have been studied extensively \cite{Causin2005,Formaggia2001, Vierendeels2005, Badia2008, Joosten2009} and various techniques to stabilize and accelerate the coupling have been proposed.
The simplest approach is a relaxation of the interface data, which improves stability at the expense of slow convergence.
Dynamically updating the relaxation factor as in Aitken relaxation \cite{Mok2001,Mok2001b,Kuttler2008} mitigates this drawback, but
it still treats all error components identically,
while it has been shown \cite{Degroote2008, Degroote2009b} that the added-mass instability is only caused by the lowest wave number components of the error between the correct solution and the one in the iterative coupling.

This realization opened the door to quasi-Newton techniques, which update the interface data communicated between the subproblems using a Newton-Raphson approach with a low-rank approximation of the Jacobian, based on input-output data from previous solver calls.
Throughout the years, many different variants have been proposed.
They differ in when the interface data is updated, i.e., after only one of the solvers or after both,
and in the technique employed to approximate the Jacobian,
such as the least-squares or the 
multi-vector approach \cite{Delaisse2021}.
An extensive overview is given in \cite{Delaisse2023}.

In literature, the efficiency of a coupling algorithm is generally assessed based on 
the required number of coupling iterations per time step,
justified by the observation that the solution of the subproblems is by far the most expensive part of the coupling scheme.
Although never stated explicitly, however, this cost measure also implies that the computational cost of each coupling iteration, and therefore each call of a solver, is constant.

This work demonstrates that this is not the case in practice.
%
Since the subproblems are typically nonlinear,
they introduce their own internal iteration loop to handle this nonlinearity, using
fixed-point or Newton iterations.
These internal \textit{subproblem iterations} have a considerable impact on the computational cost of a solver call.
As an illustrative example, imagine two solid solver calls:
the first performs ten Newton iterations, the second only one;
%
%
clearly, the second call is expected to consume much less time, although most probably not by a factor of ten.
%
%
Against this backdrop, this work proposes a new cost measure that considers not only the number of coupling iterations, but also the number of subproblem iterations in each of the two solvers, and combines them in a weighted sum.
The numerical results confirm that this new measure represents the actual wall-clock time much more accurately than looking at the coupling iterations alone.

As the impact of the \subproblemIter s on the computational cost has been overlooked so far, literature offers plenty of techniques to converge in fewer coupling steps, 
but, to the best of the authors' knowledge, lacks any studies on how to minimize the number of \subproblemIter s.
%
To start closing this gap,
this work investigates the effects of limiting the number of iterations performed for one solver call, 
demonstrating that a significant speed-up can be obtained by not converging to the final subproblem tolerance in each solver call.
This is mainly because it avoids investing time into polishing a preliminary subproblem solution that will be overwritten in the next coupling iteration, as long as the partitioned scheme has not converged yet.
On the other hand, passing back inaccurate results
brings the risk of compromising the coupling loop's stability,
as well as the quality of the input-output data used by quasi-Newton methods.

In any case,
at the end of each time step,
an accurate solution of the coupled problem requires
both subproblems to be converged up to their respective solver tolerances.
While this condition is inherently satisfied when iterating to full convergence in every solver call, limiting the number of \subproblemIter s per call requires to monitor it explicitly.
Therefore, a new convergence criterion is introduced that evaluates the convergence of the coupling scheme solely based on the solver residuals.
%
As an added benefit, this new approach avoids choosing the value of the rather non-intuitive coupling tolerance commonly used in literature
and provides a natural link between the accuracy of the coupling on the one hand and each of the subproblems on the other hand.


As different discretization approaches for solving the subproblems are common in modern engineering science,
this work discusses numerical results obtained from different partitioned FSI frameworks to broaden its scope.
While the first framework uses finite elements for both subproblems, the second employs a finite-volume method for the flow problem and finite elements for the structural problem. \\

The remainder of this work is structured as follows.
After the introduction of the fluid-structure interaction problem in \Sec{fsi}, \Sec{subproblems} treats the solution of the subproblems with finite elements or finite volumes and focuses on how both techniques manage nonlinearities.
This is important for \Sec{convergence}, which presents a new cost function, taking into account the number of coupling iterations as well as the number of \subproblemIter s. Moreover, it discusses the impact the number of \subproblemIter s per solver call has on this cost measure and
introduces the new convergence criterion.
Results are generated with both frameworks and discussed in \Sec{results} for two test cases, the lid-driven cavity and flexible tube case, before the conclusions are presented in \Sec{conclusion}.

\section{Partitioned fluid-structure interaction}
\label{sec:fsi}

For the greater part,
the research questions investigated in this work
are expected to affect
any multi-field problem
solved in a partitioned manner.
Nevertheless, 
we restrict ourselves to 
the interaction of an incompressible fluid in the domain $\flowDomain \inR{\nsd}$ and 
an elasto-dynamic solid $\solidDomain \inR{\nsd}$,
where $\nsd$ is the number of spatial dimensions.

\subsection{Incompressible flow}

The flow velocity $\fvel \ofxt$ and the fluid pressure $p^f \ofxt$
are governed by
the unsteady Navier-Stokes equations
%
%
\begin{subequations} \label{eqs:INS_Strong}
	\begin{alignat}{2} 
		\frac{\partial \fvel}{\partial t} + \fvel \cdot \boldsymbol{\nabla} \fvel  - \frac{1}{\rho^f} \boldsymbol{\nabla}  \cdot \CauchyStress^f &= \vekt{b}^f	\qquad	&& \text{in} ~\Omega^f ~~\foralltime~,\\
		\boldsymbol{\nabla}  \cdot \fvel \,&= 0 && \text{in} ~\Omega^f ~~\foralltime~,
	\end{alignat}
\end{subequations}
with the constant fluid density $\rho^f$ and the external body force $\vect{b}^f$.
%
For an incompressible Newtonian fluid,
the Cauchy stress tensor is given by $\CauchyStress^f( \fvel, p^f) = - p^f \vekt{I} + \mu^f \left(  \nabla \fvel + (\nabla \fvel)^T \right)$,
where $\mu^f$ is the dynamic viscosity.
The problem is closed by an appropriate set of boundary conditions on $\partial \Omega^f$
and a divergence-free initial velocity field.

\subsection{Elastic solid}

The solid displacement field $\vekt{d}^s(\vekt{x},t)$ is determined from 

\begin{alignat}{2} \label{equ:Elasto_Strong}
	\frac{d^2 \vekt{d}^s}{dt^2} - \frac{1}{\rho_0^s} \boldsymbol{\nabla}_0 \cdot \left( \mat{S}^s \mat{F}^T \right) &= \vekt{b}^s \qquad &&\text{in } \Omega_0^s ~~\foralltime~,
\end{alignat}
with the solid density $\rho_0^s$ and the external body force $\vect{b}^s$.
%
%
Following a total Lagrangian viewpoint, this equation of motion is formulated with respect to the
undeformed reference state $\Omega_0^s$, indicated for all affected quantities and operators by the subscript $0$.
Accordingly, the inner stress is not expressed in terms of the
Cauchy stress tensor $\CauchyStress^s$, but the second Piola-Kirchhoff stress tensor $\mat{S}^s = \det(\mat{F})\, \mat{F}^{-1} \CauchyStress^s \mat{F}^{-T}$
with the deformation gradient $\mat{F}$.
%
As constitutive equation, a Hookean or the St. Venant-Kirchhoff material model is used, resulting in a geometrically nonlinear solid problem \cite{Bathe1996}.
It is closed by an initial (zero) displacement field and suitable boundary conditions on $\partial \Omega^s$.

\subsection{Coupling conditions}

To ensure the conservation of mass, momentum, and mechanical energy over the 
shared interface  
$\Gamma^{fs} = \partial \flowDomain \cap \partial \solidDomain$,
the solution fields of the two subproblems have to satisfy
kinematic and dynamic continuity:
\begin{subequations}
\begin{align}
    \vekt{d}^f &= \vekt{d}^s && \text{on } \Gamma^{fs} ~~\foralltime ~, \label{equ:kinCont} \\
    \CauchyStress^f  \cdot \vekt{n}^f  &= - \CauchyStress^s \cdot \vekt{n}^s  && \text{on } \Gamma^{fs} ~~\foralltime ~,
\end{align}
\end{subequations}
where $\vect{d}^f$ is the fluid's displacement, while $\vekt{n}^{f}$ and $\vekt{n}^{s}$ are the interface unit normal vectors pointing outwards from the corresponding domains.
Note that \Equ{kinCont} implies the equality of velocities and accelerations too. 

\subsection{Dirichlet-Neumann partitioning}

In partitioned fluid-structure interaction simulations,
the two subproblems are addressed by two distinct solvers that are coupled in a black-box manner, i.e., solely via the exchange of interface data.
This strategy is very flexible and modular concerning the solvers,
but their communication entails some additional challenges.
On the one hand, the interface discretizations of the two subproblems in general do not match, so that transferring data requires a spatial projection.
On the other hand, 
an iterative procedure is needed to find consistent solutions of the
two interdependent subproblems within each time step.

The most common partitioned approach is the combination of a \marktext{Dirichlet-Neumann} partitioning with a Gauss-Seidel type iteration 
scheme\footnote{Although all numerical experiments of this work use a Dirichlet-Neumann partitioning, the investigated impact of the \subproblemIter s on computational cost is expected to be essentially the same for other partitionings, such as Robin-Neumann or Robin-Robin schemes \cite{Badia2008,nobile2008effective,spenkeRNQN}.}.
For every coupling iteration, it solves the flow problem with the current interface deformation and passes the interface tractions 
as a Neumann boundary condition to the solid.
The solid solver then computes the new deformation state and returns the interface displacement
to the flow solver,
where the resulting interface velocity 
poses a
Dirichlet condition.
Once this procedure has converged, the next time step is started.

The main drawback of the Dirichlet-Neumann partitioning is its sensitivity to the added-mass effect. 
%
For this work, it is sufficient to note that this instability 
is inherent to partitioned solution schemes and
increases mainly with the density ratio $\rho^f/\rho^s$, but more
detailed investigations can be found in literature \cite{Causin2005, Forster2007, VanBrummelen2009}.

\subsection{Interface quasi-Newton methods}
\label{sec:IQN}

An effective countermeasure against the added-mass effect is
to modify the interface deformation with an \marktext{interface quasi-Newton} (IQN) method
before passing it to the flow solver\footnote{Although it is a lot less common, updating the interface tractions before passing them to the solid solver is possible too \cite{spenkeRNQN}.}.
Identifying the solution of the coupled problem as a fixed-point of the coupling iteration loop,
their basic idea is to employ the Newton-like update step
\begin{align}
    \vect{d}^{k+1} = 
    \vect{\tilde{d}}^k + \Delta\vect{\tilde{d}}^k_{IQN} =
    \tilde{\vect{d}}^k - \left( \frac{\partial \Rk{k}}{\partial \tilde{\vect{d}}^k} \right)^{-1} \Rk{k} ~,
\end{align}
where $\tilde{\vect{d}}^k$ is the interface deformation computed by the solid solver 
and $\vect{d}^{k+1}$ the one sent to the flow solver, before starting the next coupling iteration $k+1$. 
$\Rk{k} \equiv \tilde{\vect{d}}^k -\vect{d}^k$ denotes the fixed-point residual of the interface deformation.

Since the exact Jacobian $\frac{\partial \Rk{k}}{\partial \tilde{\vect{d}}^k}$ is not available for black-box solvers, however,
a low-rank approximation is used instead.  
To avoid additional solver calls, this inverse Jacobian approximation $\Jac  \approx \left( \frac{\partial \Rk{k}}{\partial \tilde{\vect{d}}^k} \right)^{-1}$ is constructed from the interface deformation states that were computed in previous coupling iterations. 

For details on the concept of interface quasi-Newton methods, different variants, and implementation aspects
the authors recommend the works \cite{Delaisse2023, spenkeRNQN, spenke2020multi, Delaisse2022, Lindner2015, Degroote2009}.

\section{Solution of nonlinear subproblems}
\label{sec:subproblems}

Partitioned FSI couples two black-box solvers: one for the flow and one for the solid problem.
To increase the scope of this work, two different sets of solution techniques for the subproblems are investigated, both of which are widely used in FSI.
The first one uses finite-element (FE) methods for both the flow and solid problem, whereas the second framework relies on a finite-element method for the solid, but a finite-volume (FV) method for the flow problem.
Both setups are summarized below:

\begin{itemize}
    \item In the first simulation framework, the in-house solver XNS discretizes the flow problem by stabilized Lagrangian P1-P1 finite elements in space \cite{Pauli2017,donea2003finite} and a backward Euler method in time \cite{forti2015semi}. The ALE mesh is adapted to deforming domains via the linear elastic mesh-update method \cite{johnson1994mesh,behr2002free},
    and its velocity determined by a first-order finite difference scheme \cite{forster2006geometric}. The solid subproblem is solved by the in-house code FEAFA using isogeometric analysis \cite{hughes2005isogeometric,cottrell2009isogeometric}, a spline-based variant of finite elements, in space and a generalized-$\alpha$ scheme in time \cite{chung1993time,erlicher2002analysis}. Non-matching interface discretizations are handled by a spline-enhanced version of finite-interpolation elements \cite{hostersspline, makespline}.
    This setup will in the following be labeled \textbf{FE-FE}.
    \item The second framework uses a finite-volume method for the flow problem and a backward Euler discretization in time within the commercial solver ANSYS Fluent \cite{Fluent2019R3}.
    The mesh is structured and the discretization scheme for the convection terms of the momentum equations is second-order upwind \cite{Barth1989}.
    For the pressure equation, the second-order and the standard scheme \cite{Rhie1983} are used for the lid-driven cavity and flexible tube case, respectively, see \Sec{results}.
    The deforming fluid domain is included using the arbitrary Lagrangian-Eulerian (ALE) frame of reference and mesh deformation is based on spring-based smoothing \cite{Batina1990}. The solid problem is discretized by piecewise linear finite elements in space and a generalized-$\alpha$ scheme in time, using the Structural Mechanics Application of the Kratos Multiphysiscs code \cite{Ferrandiz2022}. The coupling between the two is performed by the in-house code CoCoNuT \cite{Delaisse2022}. The most recent code can be found in the GitHub repository \href{https://github.com/pyfsi/coconut}{pyfsi/coconut}. Data exchange on the non-matching interface is realized with radial basis mapping \cite{Lombardi2013}. This set of solution techniques will from here on be termed \textbf{FV-FE}.
\end{itemize}

Both finite-element and finite-volume methods are common discretization schemes in modern computational engineering science.
Since they are well-documented in literature, see for example \cite{Mueller2015, Versteeg1995, ferziger2002} for finite volumes and \cite{Bathe1996, zienkiewicz2005finite, zienkiewicz2005femsolids, reddy2019introduction} for finite elements, any in-depth discussion is omitted here for the sake of conciseness.

The aspect most important for this work is that both techniques transform a continuous partial differential equation (PDE) and its boundary conditions into a discrete set of algebraic equations.
In general, this system of equations is nonlinear, meaning that the system matrix $\Am$ depends on the solution $\uv$, yielding the matrix form
\begin{equation}
    \label{equ:nonlinearSystem}
    \Am(\uv) \, \uv = \bv,
\end{equation}
where $\Am \inR{n \times n}$ is the sparse coefficient matrix, $\uv \inR{n}$ holds the dependent variables of interest, 
and $n$ denotes the number of degrees of freedom (DOF).
%
%
Note that the right-hand side (RHS) vector $\bv$, containing source terms and boundary conditions,
is considered independent of $\uv$ within each subproblem\footnote{In many algorithms, part of the dependence on $\uv$ is treated explicitly resulting in a lagging contribution to the RHS $\bv$. Since this is not the focus in this work and for the sake of simplicity, it is assumed that all dependence on $\uv$ is treated implicitly, i.e., within the system matrix $\Am(\uv)$.\label{fnt:implicit}}. \Sec{convergence_criterion} will explain, however, that when coupling two subproblems as in FSI, the RHS in fact becomes a function of $\uv$ as well, i.e., $\bv = \bv(\uv)$.
For unsteady PDEs, such as \Eqs{INS_Strong} or \Equ{Elasto_Strong}, additionally a time-stepping scheme is applied, so that an algebraic system of the form \Equ{nonlinearSystem} is obtained for every time step.
It can then be solved using suitable numerical methods.


The following subsections focus on how the nonlinearity in the discrete system of equations \Equ{nonlinearSystem} is treated.
Although the system of equations can be solved using Newton or fixed-point iterations in both finite-element and finite-volume methods, this work follows the most common approach: 
Newton iterations for the finite-element method and fixed-point iterations for the finite-volume method.
The most important takeaway is that finite-element method solves \Equ{nonlinearSystem} with a few Newton iterations (typically 2 to 5), while the finite-volume method uses much cheaper fixed-point iterations to reach the solution, but requires a larger number of them (typically 50 to 200).

\subsection{Finite-element method} \label{sec:FE_Newton}

Although originally developed for solid mechanics,
finite-element methods are widely used in many fields of computational engineering,
including fluid dynamics \cite{donea2003finite,reddy2010finite}.
After transforming the PDE into its variational form,
the domain is divided into a mesh of \textit{elements}
to approximate the solution by a linear combination of the elements' basis functions.
The unknown coefficients
are then determined from the resulting algebraic system of the form \Equ{nonlinearSystem}.

In finite-element methods, it is common practice to employ Newton's method to tackle the nonlinearity of $\Am$.
Starting from an initial guess $\mathbf{u}^0$, 
in each iteration $i=1,\cdots$
the linearized system
\begin{equation}
    \label{equ:linearizedFE}
    \underbrace{ \left( \Am^{i-1} + \left.\frac{\partial \Am}{ \partial \uv} \right\vert_{i-1} \, \uv^{i-1} \right)}_{=: \mat{K}^{i-1}} ~\Delta \uv^i 
    = \underbrace{ \bv - \Am^{i-1} \, \uv^{i-1} }_{=: \resProblem^{i-1}}
    \quad \Longleftrightarrow \quad
    \mat{K}^{i-1} \, \Delta \uv^i = \resProblem^{i-1} 
\end{equation}
 is solved for the solution increment $\Delta \uv^i$, where $\mat{K}^{i-1}$ is the tangent stiffness (or system) matrix,
and
$\Am^{i-1}$ is a shorthand notation for $\Am(\uv^{i-1})$,
Further,
$\resProblem^{i-1}$ is the residual vector of the considered subproblem,
i.e., either the fluid ($p=f$) or the solid ($p=s$) problem.
Subsequently, the solution field is updated by $\uv^{i} = \uv^{i-1} + \Delta \uv^i$.
The Newton iteration is considered converged when the residual norm is lower than some tolerance $\varepsilon$, i.e., if $\frac{\normE{\resProblem^i}}{\sqrt{n}}< \varepsilon$.

Note that the computational cost of this procedure is typically dominated by the assembly of the linear system in \Equ{linearizedFE} on the one hand and its numerical solution on the other hand (in this work via a preconditioned GMRES \cite{saad1986} approach). 
Both these operations are repeated for every Newton iteration.

\subsection{Finite-volume method}
Finite-volume methods are very common for solving flow problems.
The principle of the finite-volume method is to discretize the spatial domain into \textit{finite volumes} or cells and apply the integral form of the governing equations on each of them.
In these integral forms, the volume integral of the divergence is transformed into a surface integral over its boundaries, using Gauss's theorem, so that this method is conservative by design.
Finally, this results in a large system of algebraic equations in the form of \Equ{nonlinearSystem}.
In pressure-based finite-volume solvers, pressure and velocity are either solved for together, in a so-called coupled approach (as done in this work), or sequentially, in a segregated approach. In the latter case, \Equ{nonlinearSystem} is only a symbolic notation.

The system of equations can also be solved with Newton iterations, but it is more common to linearize the nonlinear coefficient matrix $\Am$ using Picard or fixed-point iterations
\begin{equation}
    \label{equ:linearizedFV}
    \Mm^{i-1} \, (\uv^{i}-\uv^{i-1}) = \bv - \Am^{i-1} \, \uv^{i-1},
\end{equation}
where $\Mm^{i-1}=\Mm(\uv^{i-1})$ is some approximation of $\Am^{i-1}$, e.g., $\mathrm{diag}\,(\Am^{i-1})$.
It should be noted that this notation is only symbolic, as the actual solution techniques are more involved and therefore out of scope for this work, see for example the multigrid methodologies, e.g., the algebraic multigrid method (AMG) with smoother \cite{Stueben2001}.
The iteration is considered converged when the residual norm is lower than some tolerance $\varepsilon$.
In this work, the residual is calculated as the total imbalance scaled by a factor representative for the flow rate of the respective component of $\uv$ through the domain.

\section{Computational cost and convergence of partitioned algorithms} 
\label{sec:convergence}

Despite the technical differences illustrated in the previous section, both finite-element and finite-volume solvers tackle nonlinear problems by use of an internal iteration loop that repeats the assembly and (partial) solution of a linearized equation system for each iteration.
For the sake of a simpler nomenclature, in the following the term \textit{\subproblemIter s} will be used to refer to either Newton iterations or fixed-point iterations.

This section discusses the influence of these \subproblemIter s on the computational cost of partitioned schemes for fluid-structure interaction.
Moreover, it proposes a new criterion to evaluate the coupling loop's convergence solely based on the subproblem residuals.\\

\mypara For a clear notation, 
this work uses an uppercase $N$ for all iteration counts that refer to the whole simulation, like the total number of coupling iterations
$\totalCoupleIter$.
%
In contrast, a lowercase $n$ stands for the number of \subproblemIter s within a flow or solid solver call,
$n^f$ and $n^s$.
%
The \subproblemIter\ index $i$ restarts from $1$ for each solver call,
the coupling iteration index $k$ for each time step.
%
%
Additionally,
for a simpler derivation
\Sec{costMeasure} uses
the index $\kbar$ to iterate
over all $\totalCoupleIter$ coupling steps 
of the full simulation.


%

\subsection{Cost measure} \label{sec:costMeasure}

It is common practice in partitioned fluid-structure interaction simulations to treat the two solvers as black boxes
of which only the in- and output are known, but not their interior properties and functionalities.
The computational cost of a specific scheme is then quantified by the number of coupling iterations required for convergence.
Although this measure is motivated by the assumption that solving the subproblems is much more expensive than the data exchange or any other step of the coupling,
it in fact also presumes the cost of one solver call to be constant.

This may not be in line with empirical observations,
but 
following the black-box perspective in the strictest sense possible, i.e., 
if the interior workings of the solvers are completely unknown,
this is indeed the best guess.
In practice, however,
it is almost always known 
whether the subproblems are nonlinear and whether
the solvers use any standard discretization technique
like finite volumes or elements.

As discussed in \Sec{subproblems}, both FV and FE solvers 
repeat the most expensive steps, the assembly and (partial) solution of the linearized system, for every \subproblemIter.
It is therefore expected that a significant part of the cost of a solver call scales with the number of \subproblemIter s it performs.
Consequently, this work proposes to estimate the computational cost of calling a subproblem $p$ in coupling iteration $\kbar$ by
\begin{align}
    \costSolver^p(\kbar) \approx \costFix^p + {n^p}(\kbar) \cdot \costIter^p ~, \label{equ:CostSolverCall}
\end{align}
where 
$n^p(\kbar)$ is the number of \subproblemIter s run for this solver call.
The constant cost $\costFix^p$ accounts for all operations that are executed
once per solver call, 
which can for example include updating the mesh or computing the fluid loads.
In contrast, $\costIter^p$ represents all cost contributions that are incurred for each subproblem iteration,
i.e., in particular assembling and solving the linear system.
%
%
%
Note that the costs factors $\costFix^p$ and $\costIter^p$ are assumed independent from the iteration or time step.
The validity of this assumption will be confirmed by the results in \Sec{results}.

Concerning the full run time of the simulation $\timeSimulation$,
the cost measure $\costSimulation$ is obtained by summing up the cost of solver calls and data transfer over the total number of coupling iterations $\totalCoupleIter$, yielding
\begin{align}
    \timeSimulation \approx \costSimulation = \sum_{\kbar=1}^{\totalCoupleIter} \left[ \costCouple + \sum_{p=f,s}  \costSolver^p(\kbar) \right] = \totalCoupleIter \cdot \costCouple + \sum_{\kbar=1}^{\totalCoupleIter} \sum_{p=f,s}  \costSolver^p(\kbar) ~, \label{equ:CostSimulation}
\end{align}
where the cost $\costCouple$
of the data transfer and update techniques,
like relaxation or IQN, per coupling iteration
was considered 
constant. 
Note that the index $\kbar$ running over all
coupling iterations is used to avoid
introducing an additional time step sum. 

The cost measure used in this work is obtained by inserting \Equ{CostSolverCall} into \Equ{CostSimulation} and regrouping the terms based on whether they are scaling with the number of coupling iterations or the number of \subproblemIter s:
\begin{alignat}{2} \label{equ:equivalentTime}
    \costSimulation & = 
    \totalCoupleIter \cdot \costCouple ~~+~ \sumCplIter \sumProblems \big[ \costFix^p && ~+~ n^p(\kbar)  \cdot \costIter^p \big] 
    \\ \nonumber
    & = \totalCoupleIter \cdot \underbrace{\big( \costCouple  ~+~ \sumProblems \costFix^p \big)}_{=: \,\costCoupleSum} &&~+~ \sumProblems \costIter^p \underbrace{\sumCplIter n^p(\kbar)}_{=: \, \totalProblemIter^p} \\ \nonumber
    & = \totalCoupleIter \cdot \costCoupleSum ~+~ \sumProblems \totalProblemIter^p \cdot \costIter^p       
\end{alignat}
The newly introduced variables $\totalProblemIter^f$, $\totalProblemIter^s$, and $\costCoupleSum=\costCouple+\costFix^f+\costFix^s$ represent the total number of \subproblemIter s of problem $f$ or $s$ 
and all costs occurring once per coupling iteration, respectively.

The cost factors $\costIter^p$ and $\costCoupleSum$
are prescribed by the simulation framework, computer architecture, problem sizes, etc. and therefore considered
as given constants,
so that the cost measure boils down to a weighted sum of the iteration counts $\totalCoupleIter$, $\totalProblemIter^f$, and $\totalProblemIter^s$.
%
An efficient partitioned scheme therefore not only has to reduce the total number of coupling iterations, but in particular also the
total number of \subproblemIter s. \\

\mypara In partitioned FSI, the mesh update of the flow problem is sometimes interpreted as a third problem besides fluid and structure.
While this viewpoint could easily be integrated into \Equ{CostSimulation} by adding the mesh update to the sum over the subproblems,
it will be included in the flow solver's cost within this work, for the sake of a cleaner notation.
This is equivalent as long as the number of mesh updates scales with $\totalCoupleIter$ (or $\totalProblemIter^f$).

\subsection{Iterations per solver call} \label{sec:IterationsPerCall}

A straightforward and simple way to influence the total number of \subproblemIter s $\totalProblemIter^f$ and $\totalProblemIter^s$ performed throughout the whole simulation
is to limit the number of \subproblemIter s per solver call.
More precisely, each solver call may only perform up to $\iterPerCall^f$ or $\iterPerCall^s$  \subproblemIter s,
rather than always iterating until full convergence is reached.

This naturally raises the central research question investigated in this work:
how do the maximum numbers of \subproblemIter s per solver call, $\iterPerCall^f$ and $\iterPerCall^s$,
influence the total computational cost of a partitioned algorithm?
%
The question might sound simple at first, but the relation between the \subproblemIter s per solver call and the computational cost is non-trivial:
\begin{itemize}
    \item If fewer \subproblemIter s are performed per solver call,
    the solvers exchange data more frequently, so that
    the boundary conditions at the FSI interface and with it 
    the RHS $\bv$ stay up to date,
    improving the individual quality of each \subproblemIter.
    Aside from the higher communication cost,
    however, this also brings the risk of feeding back inaccurate data into the coupling loop and, if applicable, in the quasi-Newton Jacobian approximation, 
    which may result in slower convergence, i.e., a higher $\totalCoupleIter$, or even divergence.
    
    \item Running more \subproblemIter s per solver call, on the other hand,
    increases the accuracy of the exchanged data fields
    and the IQN method's input-output pairs,
    typically resulting in a reduced number of coupling steps $\totalCoupleIter$.
    On the downside, computational time is misspend on polishing a solution for which the boundary conditions are still incorrect, so that it will be overwritten in the next coupling step anyway.
\end{itemize}
Good choices of the \subproblemIter s per flow and solid solver call, $\iterPerCall^f$ and $\iterPerCall^s$,
should balance these opposing trends to reduce the overall computational cost. \\


\subsection{Convergence criterion}
\label{sec:convergence_criterion}

It is common in partitioned FSI simulations to determine the convergence of the coupling loop by comparing a norm of the fixed-point residual $\Rk{k}$
to some tolerance $\varepsilon^c$, see for example
\cite{Kuttler2008,Degroote2013b,scheufele2017robust,Schussnig_2022}.
Although 
both \subproblemIter\ loops have to converge as well
for accurate results,
%
%
these conditions are rarely accounted for explicitly,
as they are inherently fulfilled in case all solver calls
iterate to full convergence.
Limiting the \subproblemIter s per solver call
to $\iterPerCall^f$ or $\iterPerCall^s$,
however, requires a method to make sure that both subproblems satisfy their residual tolerances $\varepsilon^f$ and $\varepsilon^s$
before going to the next time step.
Unfortunately, the influence of the fixed-point tolerance $\varepsilon^c$ is rarely discussed in literature,
let alone its interplay with the subproblem residual tolerances $\varepsilon^f$ and $\varepsilon^s$.

To tackle this issue, this work proposes a novel
convergence criterion for partitioned algorithms that does not
introduce any coupling tolerance $\varepsilon^c$.
Instead, it relies solely on the subproblem residuals.
As an added benefit, the new criterion allows for a fair comparison of convergence rates,
independent from the number of \subproblemIter s per solver call.
The remainder of this section derives the proposed convergence criterion.\\



As explained in \Sec{subproblems}, both finite-element and finite-volume solvers
iteratively solve a nonlinear problem of the form $\Am(\uv)~\uv = \bv$. 
In other words, while they account for a nonlinearity of the system matrix $\Am$, the right-hand side vector $\bv$ is assumed to stay unchanged throughout the \subproblemIter s.
In a coupled FSI problem, however, the right-hand sides of both subproblems in fact also depend on $\uv$, even if only in an implicit manner.
For example, the traction and pressure forces exerted by the fluid onto the solid change with the deformation state.
Analogically, updating the flow field alters the solid deformation and with it the fluid's boundary position and velocity.
In practice, 
this implicit nonlinearity of the RHS is impossible to explicitly account for when coupling two black-box
solvers in a partitioned algorithm.

An interesting consequence of this realization is that only the first \subproblemIter\ of a solver call uses the correct RHS vector $\bv$, because
in the first iteration, both the system matrix $\Am$ and RHS $\bv$ refer to the current solution $\uv^0$, i.e., the initial value for the subproblem solve.
The subsequent \subproblemIter s, on the other hand, inevitably lack any contribution of the other problem to $\bv(\uv)$
that would
follow from the change of $\uv$ within the current solver call, as only the system matrices are updated to the new solution $\uv^i$, while $\bv=\bv(\uv^0)$ stays unchanged\footnote{As mentioned in \Fnt{implicit} (page~\pageref{fnt:implicit}), it is assumed that all dependence on $\uv$ inherent to the subproblem is treated implicitly, i.e., within the system matrix $\Am(\uv)$.}.
%
The effect is sketched in \Tab{NonlinearRHS}. 

\begin{table}[h]
	\caption{As the dependency of the RHS vector $\bv$ on the solution $\uv$ is resulting from the other subproblem and hence missing during the \subproblemIter s in a partitioned FSI simulation, only the first {\subproblemIter} of each solver call uses the correct RHS. For $i>1$, the RHS vector is defective as $\bv(\uv^0) \neq \bv(\uv^{i-1})$.} \label{tab:NonlinearRHS}
    \begin{center}
        \smallIfElsevier 
        \begin{tabular}{c c c}
        Iteration & Fixed-point iteration & Newton iteration \\
        \hline 
        $i=1$   
        & $~\,\,\Mm(\uv^0) \, (\uv - \uv^{0}) ~\,\, = \bv(\uv^{0}) - ~\,\, \Am(\uv^{0}) \, \uv^{0} ~\,\,$
        & $~\,\,\mat{K}(\uv^0) \, \Delta \uv = \bv(\uv^0) - ~\,\, \Am(\uv^0) \, \uv^0 ~\,\,$
        \\
        $i>1$   
        &  $\Mm(\uv^{i-1}) \, (\uv - \uv^{i-1}) = \bv(\uv^{0}) - \Am(\uv^{i-1}) \,\uv^{i-1}$
        &  $\mat{K}(\uv^{i-1}) \, \Delta \uv = \bv(\uv^0) - \Am(\uv^{i-1}) \,\uv^{i-1}$
    	\end{tabular}
    \end{center}
%
%
\end{table}

As a result, only the residual of a solver call's first \subproblemIter, 
$\vekt{r}_p^0 = 
\bv (\uv^0) - \Am (\uv^0) \uv^0$, allows to draw conclusions about the convergence of the coupling loop, because it quantifies to which extent the subproblem at hand already balances the current coupling data, i.e.,
the coupling data that resulted from the current subproblem solution.
%
For example, the solid solver's first residual is the difference between the Cauchy stresses and the
external loads that the flow solver determined with the current deformation state as boundary condition.
The \subproblemIter s $i=2,3,\cdots$, in contrast, merely attempt to find the converged solution for a defective right-hand side vector $\bv = \bv(\uv^0) \neq \bv(\uv^{i-1})$.

This motivates the idea to trigger the convergence of the coupled problem via the subproblem residuals of the first \subproblemIter:
the coupling iteration is considered converged, if for all subproblems the residual of the first \subproblemIter\ satisfies the respective tolerance $\varepsilon^f$ or $\varepsilon^s$.
The different iteration loops and the new convergence criteria are illustrated by \Fig{IterationLoops}.

\begin{figure}[ht]
	\centering
	\resizebox{0.99\textwidth}{!}{
		\tikzset{every picture/.style={line width=0.75pt}} 

\begin{tikzpicture}[x=0.75pt,y=0.75pt,yscale=-1,xscale=1]

\draw [line width=1.5]    (310,128) -- (310,68) -- (800,68) -- (800,124) ;
\draw [shift={(800,128)}, rotate = 270] [fill={rgb, 255:red, 0; green, 0; blue, 0 }  ][line width=0.08]  [draw opacity=0] (11.61,-5.58) -- (0,0) -- (11.61,5.58) -- cycle    ;
\draw [line width=1.5]    (801,443) -- (801,503) -- (675,503) ;
\draw [shift={(671,503)}, rotate = 360] [fill={rgb, 255:red, 0; green, 0; blue, 0 }  ][line width=0.08]  [draw opacity=0] (11.61,-5.58) -- (0,0) -- (11.61,5.58) -- cycle    ;
\draw [line width=1.5]    (310,480) -- (310,444) ;
\draw [shift={(310,440)}, rotate = 90] [fill={rgb, 255:red, 0; green, 0; blue, 0 }  ][line width=0.08]  [draw opacity=0] (11.61,-5.58) -- (0,0) -- (11.61,5.58) -- cycle    ;
\draw [line width=1.5]    (560,523) -- (560,559) ;
\draw [shift={(560,563)}, rotate = 270] [fill={rgb, 255:red, 0; green, 0; blue, 0 }  ][line width=0.08]  [draw opacity=0] (11.61,-5.58) -- (0,0) -- (11.61,5.58) -- cycle    ;
\draw  [dash pattern={on 1.69pt off 2.76pt}][line width=1.5]  (281,480) -- (341,480) -- (341,523) -- (281,523) -- cycle ;
\draw [line width=1.5]    (450,502) -- (344,502) ;
\draw [shift={(340,502)}, rotate = 360] [fill={rgb, 255:red, 0; green, 0; blue, 0 }  ][line width=0.08]  [draw opacity=0] (11.61,-5.58) -- (0,0) -- (11.61,5.58) -- cycle    ;


\draw  [line width=1.5]   (118,128) -- (502,128) -- (502,442) -- (118,442) -- cycle  ;
\draw (310,285) node   [align=left] {\begin{minipage}[lt]{258.4pt}\setlength\topsep{0pt}
\large 
\setstretch{1.15}
\begin{center}
\textbf{Flow solver}
\end{center}
\vspace{0.1cm}
$\text{converged}^f=false$ \tabto{5cm} \textcolor{gray}{// initialization}\\
\textbf{do} $\displaystyle i=1,...,n_{max}^{f}$     \tabto{5cm} \textcolor{gray}{// subproblem iteration}\\
\tabto{0.5cm} $\displaystyle \mathbf{A}^{i-1} = \mathbf{A}\left(\mathbf{u}^{i-1}\right)$  \tabto{5cm} \textcolor{gray}{// assemble system}\\	
\tabto{0.5cm} $\displaystyle \resFlow^{i-1} \,\,=\mathbf{b} -\mathbf{A}^{i-1} \ \mathbf{u}^{i-1}$ \tabto{5cm} \textcolor{gray}{// evaluate residual}\\	
\tabto{0.5cm} $\displaystyle \mathbf{u}^i \,= \mathbf{A}^{i-1} \, \backslash \, \mathbf{b}$ \tabto{5cm} \textcolor{gray}{// update solution}\\	
\tabto{0.5cm} \textbf{if} $\displaystyle \norm{\resFlow^{i-1}} < \varepsilon ^{f}$ \textbf{then} \tabto{5cm} \textcolor{gray}{// check convergence}\\
\tabto{1.0cm} \textbf{if} $i == 1$ \textbf{then} \tabto{5cm} \textcolor{gray}{// check if 1st residual} \\  	
\tabto{1.5cm} $\text{converged}^f = true$  \tabto{5cm} \textcolor{gray}{// set convergence flag}\\  
\tabto{1.0cm} \textbf{end if}       \\
\tabto{1.0cm} break                      \tabto{5cm}  \textcolor{gray}{// leave iteration loop}\\
\tabto{0.5cm} \textbf{end if}\\
\textbf{end do}\\                       
$\displaystyle \uv^{0} =\uv^{i}$    \tabto{5cm} \textcolor{gray}{// prepare next call}
\end{minipage}};
\draw  [line width=1.5]   (608,128) -- (992,128) -- (992,442) -- (608,442) -- cycle  ;
\draw (800,285) node   [align=left] {\begin{minipage}[lt]{258.4pt}\setlength\topsep{0pt}
\large 
\setstretch{1.15}
\begin{center}
\textbf{Solid solver}
\end{center}
\vspace{0.1cm}
$\text{converged}^s=false$ \tabto{5cm} \textcolor{gray}{// initialization}\\
\textbf{do} $\displaystyle i=1,...,n_{max}^{s}$     \tabto{5cm} \textcolor{gray}{// subproblem iteration}\\
\tabto{0.5cm} $\displaystyle \mathbf{A}^{i-1} = \mathbf{A}\left(\mathbf{u}^{i-1}\right)$  \tabto{5cm} \textcolor{gray}{// assemble system}\\	
\tabto{0.5cm} $\displaystyle \resSolid^{i-1} \,\,=\mathbf{b} -\mathbf{A}^{i-1} \ \mathbf{u}^{i-1}$ \tabto{5cm} \textcolor{gray}{// evaluate residual}\\	
\tabto{0.5cm} $\displaystyle \mathbf{u}^i \,= \mathbf{A}^{i-1} \, \backslash \, \mathbf{b}$ \tabto{5cm} \textcolor{gray}{// update solution}\\	
\tabto{0.5cm} \textbf{if} $\displaystyle \norm{\resSolid^{i-1}} < \varepsilon ^{s}$ \textbf{then} \tabto{5cm} \textcolor{gray}{// check convergence}\\
\tabto{1.0cm} \textbf{if} $i == 1$ \textbf{then} \tabto{5cm} \textcolor{gray}{// check if 1st residual} \\  	
\tabto{1.5cm} $\text{converged}^s = true$  \tabto{5cm} \textcolor{gray}{// set convergence flag}\\  
\tabto{1.0cm} \textbf{end if}       \\
\tabto{1.0cm} break                      \tabto{5cm}  \textcolor{gray}{// leave iteration loop}\\
\tabto{0.5cm} \textbf{end if}\\
\textbf{end do}\\                       
$\displaystyle \uv^{0} =\uv^{i}$    \tabto{5cm} \textcolor{gray}{// prepare next call}
\end{minipage}};

\draw (509,48) node [anchor=north west][inner sep=0.75pt]   [align=left] {\large {Coupling loop}};

\draw  [dash pattern={on 1.69pt off 2.76pt}][line width=1.5]   (448,480) -- (672,480) -- (672,524) -- (448,524) -- cycle  ;
\draw (560,502) node   [align=left] {\begin{minipage}[lt]{149.6pt}\setlength\topsep{0pt}
\large
\begin{center}
\textbf{if}
$\text{converged}^f \textbf{ and } \text{converged}^s$ 
\end{center}

\end{minipage}};
\draw (509,568) node [anchor=north west][inner sep=0.75pt]   
[align=left] {\large{Next time step}};
\draw (570,530) node [anchor=north west][inner sep=0.75pt]   [align=left] {\textbf{then}};
\draw (415,483) node [anchor=north west][inner sep=0.75pt]   [align=left] {\textbf{else}};

\draw (297,494) node [anchor=north west][inner sep=0.75pt]   [align=left] {\large{IQN}};

\end{tikzpicture}
	}
	\caption{Illustration of the different iteration loops considered in this work. While there is only one coupling loop managing the solver calls,
	each subproblem has its own internal iterations. 
 Note that the pseudo code $\mathbf{u}^i \,= \mathbf{A}^{i-1} \, \backslash \, \mathbf{b}$ 
 represents solving either \Equ{linearizedFE} for a Newton iteration
 or \Equ{linearizedFV} for a fixed-point iteration. Similarly, 
 $\norm{\resFlow^{i-1}} < \varepsilon ^{f}$ is merely a symbolic notation
 for the convergence checks discussed in \Sec{subproblems}. 
 Lastly, the indicated quasi-Newton update is of course optional.
 } \label{fig:IterationLoops} 
\end{figure}

To integrate the new criterion into an existing black-box FSI framework,
it is typically simplest to check after each solver call whether (1) the subproblem has converged
and (2) only one \subproblemIter\ was run. Only if both conditions are fulfilled for flow and solid solver in the current coupling iteration,
the time step has converged. \\


\mypara Note that in \Fig{IterationLoops} the solution of the subproblem is updated after the calculation of the residual, but before convergence is checked.
As a consequence, one more update is performed even though the \subproblemIter\ converged.
It would also be possible to exit the \subproblemIter\ upon convergence.
However, this is often impossible for black-box solvers and, 
moreover, the accuracy of the coupled solution would then be prescribed by
the subproblem convergence tolerance that is satisfied first, rather than the one satisfied last. \\

\mypara The interpretation of the first subproblem residual within a solver call as a measure for the coupling loop's
convergence is in fact only exact if the exchanged data fields are not modified in between the solvers.
%
%
Although this modification is common in practice
and for example part of relaxation or IQN methods (\Sec{IQN}), 
the update increment is proportional to the fixed-point residual $\Rk{k}$ in these approaches
and therefore
vanishes upon convergence, 
as will be confirmed in \Sec{results}.
Consequently, it is safe to neglect the effect of an update step within the proposed convergence criterion. \label{Remark:Increment}

\section{Results}
\label{sec:results}

This section discusses the application of the developed concepts to numerical test cases.
To investigate the influence on
the computational cost,
each test case is run multiple times with varying maximum
numbers of \subproblemIter s per solver call.

\subsection{Setup of cases}

Two different cases are studied, both well-known benchmarks for FSI.
They are simulated in the FE-FE as well as the FV-FE framework.

\subsubsection{Lid-driven cavity case}
\begin{figure}[ht]
    \centering
    \begin{subfigure}[t]{.49\textwidth}
        \input{Figure_CavityGeometry}
    	\caption{Sketch of the geometry and boundary conditions.}
    	\label{fig:CavityGeometry}
    \end{subfigure}
    \begin{subfigure}[t]{.49\textwidth}
        \centering
        \includegraphics[width=0.65\textwidth,trim={0 0 0 0}]{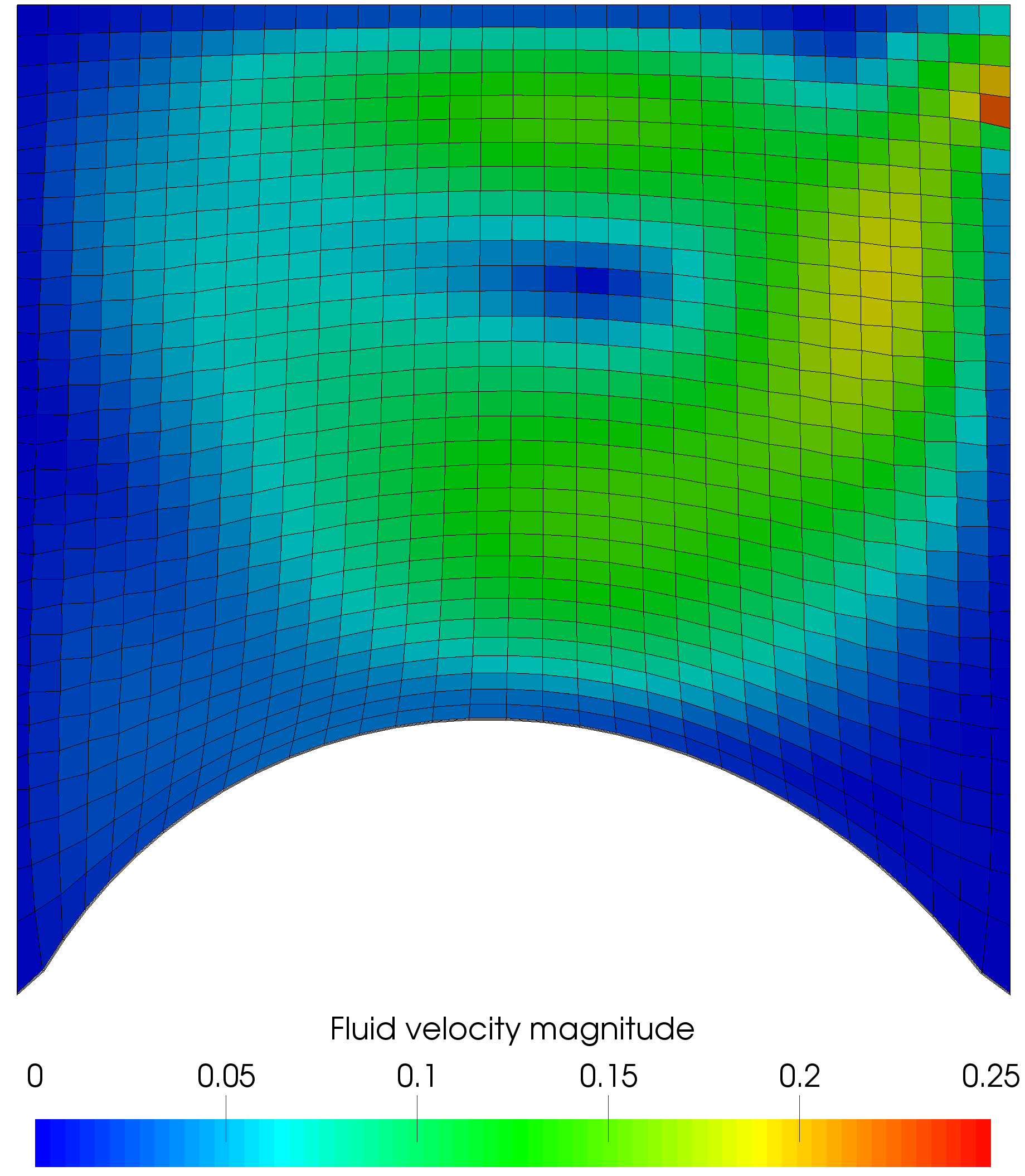}
    	\caption{Deformation of the flexible bottom and corresponding velocity field in \qty{}{\metre\per\second} at $t=\qty{40}{\second}$ simulated in the FV-FE framework.}
    	\label{fig:CavityIllustration}
    \end{subfigure}
	\caption{Visualisations of the lid-driven cavity case.}
	\label{fig:CavityCase}
\end{figure}

The first case considers a two-dimensional (2D) cavity with flexible bottom, excited by an oscillatory incoming flow and lid motion \cite{Mok2001a,Valdes2007}. 
The geometry of the cavity and the flexible bottom is detailed in \Fig{CavityGeometry} and \Tab{TubeCavityParameters} contains the structural density $\rho^s$, modulus of elasticity $E^s$, and Poisson ratio $\nu^s$ of the lid, as well as the density $\rho^f$ and dynamic viscosity $\mu^f$ of the liquid.
Both ends of the flexible bottom are fixed.
Further, on the top boundary of the cavity, a horizontal velocity $\bar{v}(t)$ is applied, which is expressed (in \qty{}{\metre\per\second}) as a function of time $t$ by
\begin{equation}
    \label{equ:CavityVelocity}
    \bar{v}=1-\cos\left(\frac{2\pi t}{5}\right).
\end{equation}
On the upper part of the left boundary, the same velocity is applied, but increasing linearly from zero to $\bar{v}(t)$ at the top, as indicated in \Fig{CavityGeometry}.
On the corresponding region of the right boundary, a pressure of \qty{0}{\pascal} is imposed.
The simulation runs 700 time steps of \qty{0.1}{\second},
adding up to the total time \qty{70}{\second}.
\Fig{CavityIllustration} illustrates the deformation of the flexible bottom.

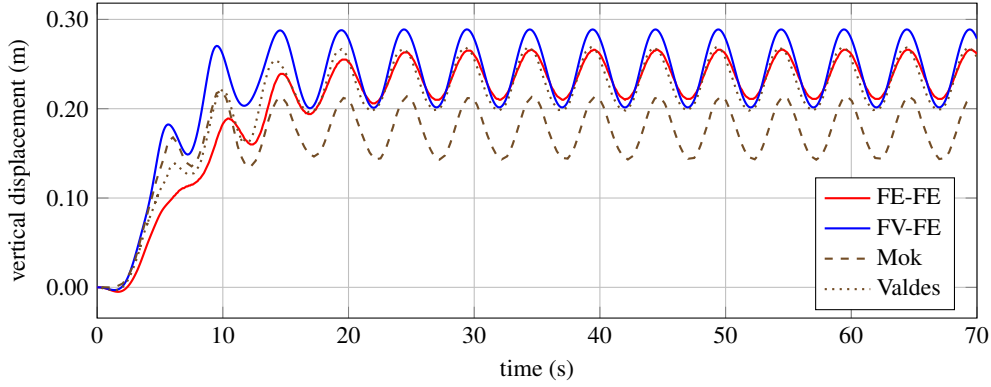
\begin{figure}
    \centering
    \begin{tikzpicture}
    \smallIfElsevier 
    \begin{axis}[
        height=0.35\textwidth,
        width=0.8\textwidth,
        xmin=0,
        xmax=70,
        grid=both,
        xlabel=time ($\qty{}{s}$),
        ylabel=vertical displacement ($\qty{}{m}$),
        /pgf/number format/.cd,
        fixed,
        y tick label style={/pgf/number format/fixed zerofill},
        legend style={
            at={(0.975,0.025)},
            anchor=south east,
            legend cell align=left,
            },
        ]
        \addplot [thick,red] table[x=time,y=uy_FE,col sep=comma]{Data/cavity_displacement.csv};
        \addlegendentry{FE-FE}
        \addplot [thick,blue] table[x=time,y=uy_FV,col sep=comma]{Data/cavity_displacement.csv};
        \addlegendentry{FV-FE}
        \addplot [thick,brown!60!black,dashed] table[x=time_mok,y=uy_mok,col sep=comma]{Data/cavity_displacement.csv};
        \addlegendentry{Mok}
        \addplot [thick,brown!60!black,dotted] table[x=time_valdes,y=uy_valdes,col sep=comma]{Data/cavity_displacement.csv};
        \addlegendentry{Valdes}
    \end{axis}
\end{tikzpicture}
	\caption{Vertical displacement of the top center point of the flexible bottom wall as a function of time. The results of both frameworks are compared with each other as well as the results obtained by Mok \cite{Mok2001a} and Valdes \cite{Valdes2007}.}
	\label{fig:CavityComparison}
\end{figure}

The results obtained with the two frameworks are compared with each other as well as with data from literature in \Fig{CavityComparison}, showing that the period of the oscillation matches well.
Although there is some variation in the amplitude, the differences between the two frameworks are in line with the ones observed in literature.

\subsubsection{Flexible tube case}
The second case is the simulation of a pressure pulse travelling through a flexible tube \cite{Degroote2009}.
In contrast to the lid-driven cavity case, this simulation is performed in three dimensions (3D).
The tube has a length of \qty{0.05}{\metre}, a radius of \qty{0.005}{\metre}, and a wall thickness of \qty{0.001}{\metre}.
The material parameters are given in \Tab{TubeCavityParameters}.

\begin{table}
	\caption{Parameter values for the lid-driven cavity and flexible tube case.}
	\label{tab:TubeCavityParameters}
	\begin{center}
        \smallIfElsevier 
		\begin{tabular}{l|ll|ll}
			& Flow parameter & Value & Solid parameter & Value \\
            \hline
            & $\rho^f$ & \qty{1}{\kilogram\per\metre\cubed}  & $\rho_s$ & \qty{500}{\kilogram\per\metre\cubed} \\
			\textbf{Lid-driven cavity} & $\mu^f$ & \qty{0.01}{\pascal\second} & $E^s$ & \qty{250}{\newton\per\metre\squared} \\
			& &                                              & $\nu^s$ & 0.0 \\
			\hline
			& $\rho^f$ & \qty{1000}{\kilogram\per\metre\cubed} & $\rho^s$ & \qty{1200}{\kilogram\per\metre\cubed} \\
			\textbf{Flexible tube} & $\mu^f$ & \qty{0.003}{\pascal\second} & $E^s$ & \qty{300000}{\newton\per\metre\squared} \\
			& &                                                & $\nu^s$ & 0.3 
		\end{tabular}
	\end{center}
\end{table}

\begin{figure}[ht]
	\centering
    \includegraphics[width=0.9\textwidth,trim={200 120 200 450},clip]{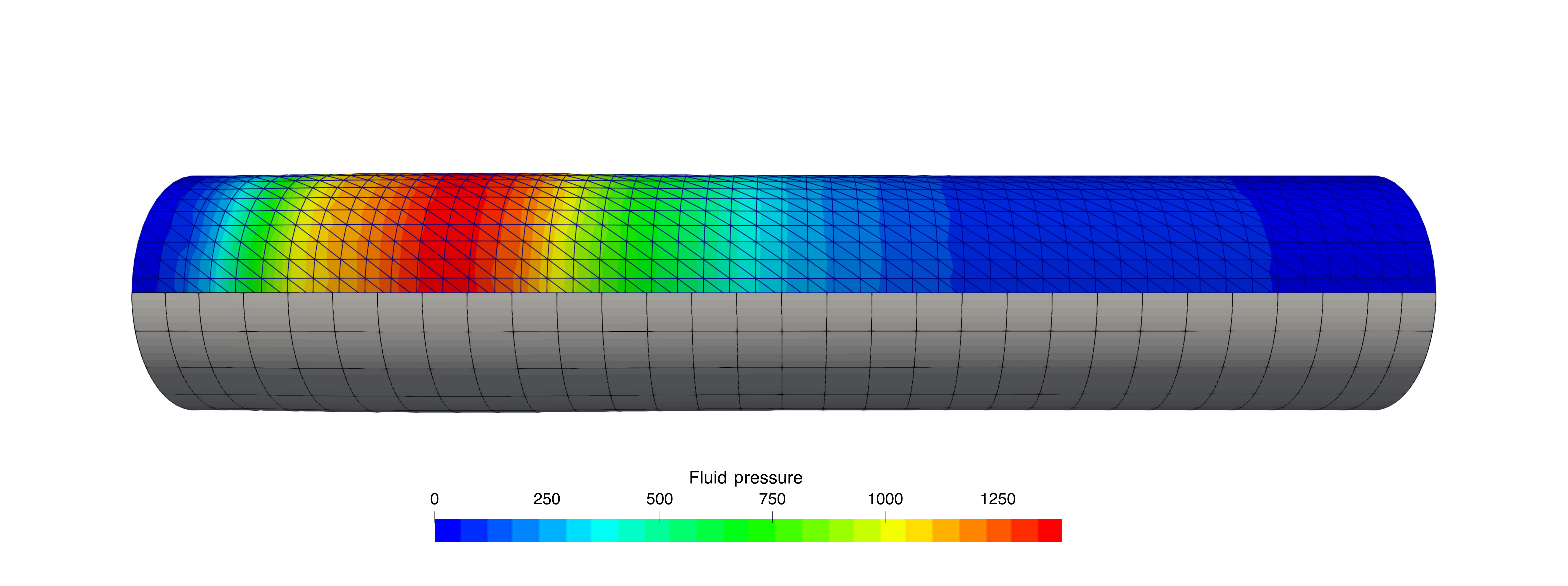}
	\caption{Flexible tube at $t=\qty{0.004}{\second}$ simulated in the FE-FE framework. The upper half shows the fluid pressure field in \qty{}{\pascal} and the lower part illustrates the IGA shell structure.}
	\label{fig:TubeIllustration}
\end{figure}

Both ends of the tube wall are clamped.
At the inlet, a pressure value of \qty{1333.2}{\pascal} is applied during the first \qty{0.003}{\second}, thereafter the pressure is \qty{0}{\pascal}.
At the outlet, the pressure is fixed to \qty{0}{\pascal}.
The total simulated time is \qty{0.01}{\second}, divided into 100 time steps of \qty{0.0001}{\second}.
\Fig{TubeIllustration} illustrates the pressure pulse travelling through the tube.

\begin{figure}[ht]
        \centering
	\begin{tikzpicture}
        \smallIfElsevier 
		\begin{axis}
			[
			height=0.35\textwidth,
			width=0.55\textwidth,
			xmin=0,
			xmax=0.05,
			grid=both,
			xlabel=axial coordinate ($\qty{}{m}$),
			ylabel=radial displacement ($\qty{}{m}$),
            x tick label style={
                    /pgf/number format/fixed,
            },
            scaled x ticks=false,
			legend style={
				anchor=north east,
				legend cell align=left,
			},
			]
			\addplot [thick,red] table[x=x_FE,y=dr_FE,col sep=comma]{Data/tube_displacement.csv};
			\addlegendentry{FE-FE}	
			\addplot [thick,blue] table[x=x_FV,y=dr_FV,col sep=comma]{Data/tube_displacement.csv};
			\addlegendentry{FV-FE}
		\end{axis}
	\end{tikzpicture}
	\caption{Radial displacement plotted over the axial coordinate at $t=\qty{0.004}{\second}$ for both solver frameworks.}
	\label{fig:TubeComparison}
\end{figure}
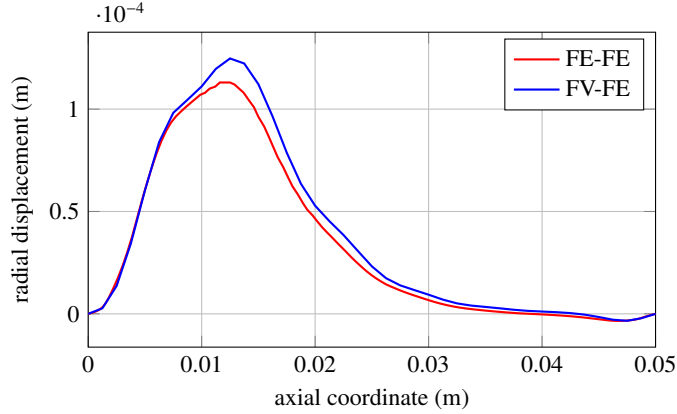

As the purpose of this study is to investigate similar effects in both software frameworks rather than a quantitative comparison of the two,
matching results are not strictly required.
Nevertheless,
\Fig{TubeComparison} compares the data obtained for both frameworks.
Despite small quantitative differences, the two graphs are very similar.

\subsubsection{Settings of coupling and solvers}

The FSI problems detailed above are solved with the quasi-Newton algorithm IQN-ILS \cite{Degroote2009}.
For each problem, the reuse parameter $q$ is detailed in \Tab{CouplingSummary}.
Before constructing the approximate Jacobian, the linear dependencies in the stored residual data are removed by filtering \cite{Degroote2013b}. 
\Tab{CouplingSummary} also contains the used filtering tolerances $\varepsilon^{fil}$, as well as the subproblem tolerances for the flow and solid solvers, $\varepsilon^f$ and $\varepsilon^s$, confer the discussion in \Sec{subproblems}.
\Tab{DiscretizationSummary} shows that similar edge lengths have been aspired for both frameworks.

\begin{table}
	\caption{Summary of the coupling and solver settings.}
	\label{tab:CouplingSummary}
	\begin{center}
        \smallIfElsevier 
        \begin{tabular}{l|*{5}{l}}
			& $q$ & $\varepsilon^{fil}$ & $\varepsilon^f$ & $\varepsilon^s$ \\
            \hline
            \textbf{Lid-driven cavity} & & & & & \\
			$\quad$ FE-FE     & 3 & $5\E{-09}$ & $1\E{-09}$ & $1\E{-09}$ \\
			$\quad$ FV-FE     & 3 & $1\E{-09}$ & $1\E{-06}$ & $1\E{-09}$ \\
			\textbf{Flexible tube}     & & & & & \\
			$\quad$ FE-FE     & 5 & $1\E{-12}$ & $1\E{-09}$ & $1\E{-09}$ \\
			$\quad$ FV-FE     & 5 & $1\E{-12}$ & $1\E{-07}$ & $1\E{-09}$ \\
		\end{tabular}
	\end{center}
\end{table}

\begin{table}
	\caption{Summary of discretization.}
	\label{tab:DiscretizationSummary}
	\begin{center}
        \smallIfElsevier 
		\begin{tabular}{l|ccc|ccc}
			& \multicolumn{3}{c|}{Flow solver} & \multicolumn{3}{c}{Solid solver} \\
			& Cells/Elements & Nodes & DOFs & Elements & Control points/Nodes & DOFs \\
            \hline
            \textbf{Lid-driven cavity} & & & & & & \\
			$\quad$ FE-FE     & 2 048 & 1 089 & 3 267 & 31 & 99 & 198 \\
			$\quad$ FV-FE     & 1 024 & 1 089 & 3 072 & 64 & 99 & 198 \\
			\textbf{Flexible tube}     & & & & & & \\
			$\quad$ FE-FE     & 35 986 & 7 195 & 28 780 & 600 & 800 & 4 800 \\
			$\quad$ FV-FE     & 7 200  & 7 913 & 28 800 & 780 & 800 & 4 800 \\
		\end{tabular}
	\end{center}
\end{table}

\subsection{Cost factors and regression}
\label{sec:regression}
As explained in \Sec{convergence} the cost of a partitioned simulation is not determined by the number of coupling iterations alone, 
but instead also depends on the number of {\subproblemIter s}.
That is why the cost measure in \Equ{equivalentTime} was introduced
as a weighted sum of the iteration counts $\totalCoupleIter$, $\totalProblemIter^f$, and $\totalProblemIter^s$.
The weights are given by a set of cost coefficients accounting for the data transfer as well as both a fixed and an iteration-dependent
contribution of each solver call.
While the values of these five cost factors are impossible to determine in general,
they can be approximated for a specific case and framework using linear regression.

To do so,
the total run time of each simulation 
is split into three parts: the time spent in the flow solver or flow time $\timing^f$, the time spent in the solid solver or solid time $\timing^s$, and finally the remainder,
which will be called coupling time $\timing^c$.
For each parameter study, i.e., a specific case and framework combination,
the regression then uses the timings and iteration counts of all conducted runs as data set.

In the following, the calculations for the flexible tube case with FV-FE are used as example for the regression. The results themselves will be described in more detail in \Sec{parameterStudy}.

\begin{figure}[ht]
    \centering
    \begin{subfigure}{.49\textwidth}
    	\includegraphics[width=\linewidth,trim={60 0 30 35},clip]{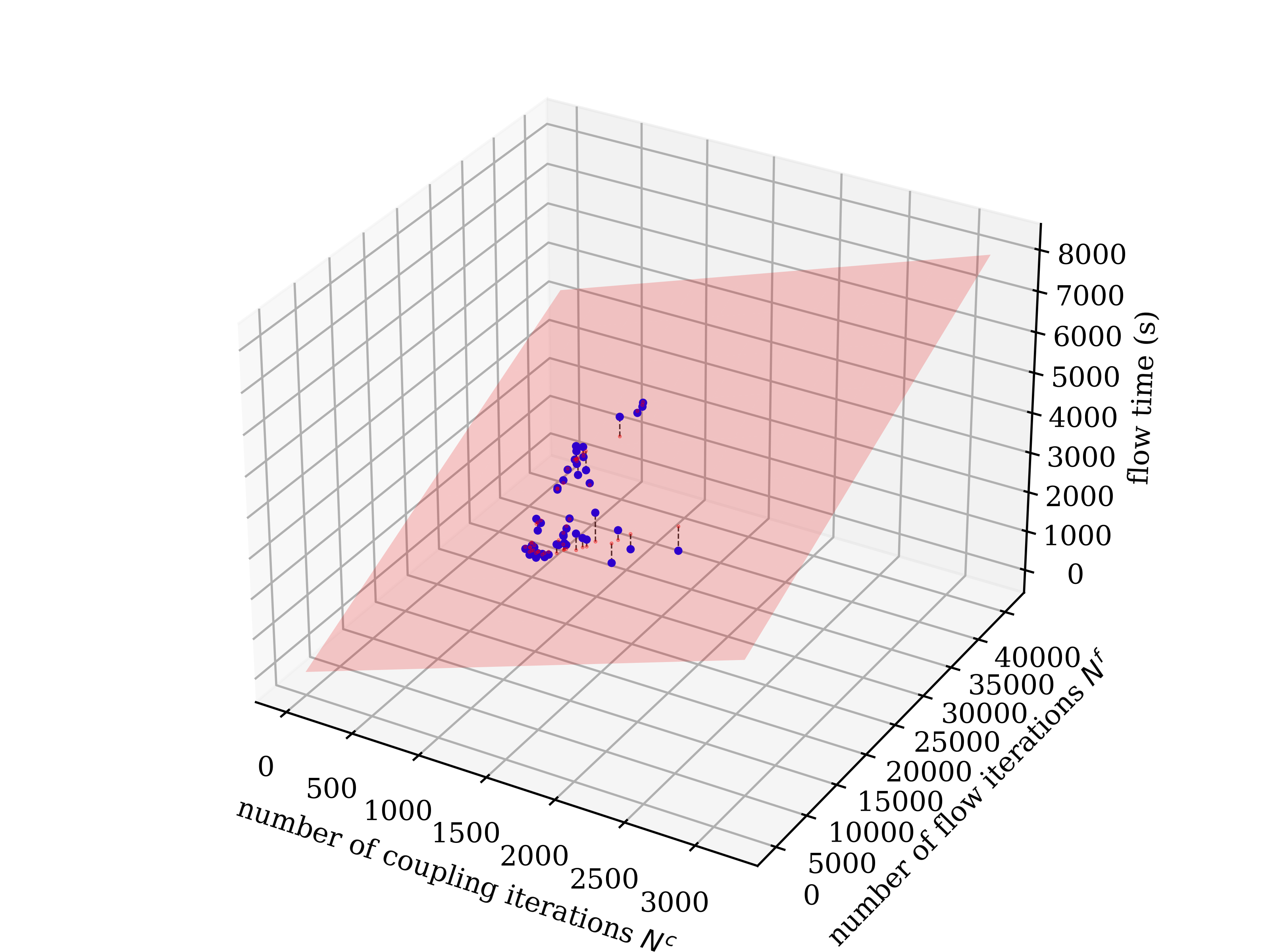}
        \caption{Plane obtained with multivariate regression.}
    	\label{fig:RegressionTubeFVFlow3D}
    \end{subfigure}
    \begin{subfigure}{.49\textwidth}
        \centering
        \begin{tikzpicture}
    \smallIfElsevier 
    \begin{axis}[
        width=1.0\textwidth,
        height=0.8\textwidth,
        grid=both,
        xmin=0,
        xmax=30000,
        ymin=0,
        ymax=5000,
        xlabel=number of flow iterations $\totalProblemIter^f$,
        ylabel=flow time ($\qty{}{\second}$),
        /pgf/number format/.cd,
        fixed,
        1000 sep={\ },
        legend style={
            at={(0.975,0.025)},
            anchor=south east,
            legend cell align=left,
            },
        ]
        \addplot [only marks,thick,blue,mark=*] table[x=it,y=time,col sep=comma]{Data/regression_flow.csv};
        \addlegendentry{Actual time $\timing^f$}
        \addplot [only marks,thick,red,mark=square*] table[x=it,y=eq_time,col sep=comma]{Data/regression_flow.csv};
        \addlegendentry{Fitted time $\totalProblemIter^c \cdot \costFix + \totalProblemIter^f \cdot \costIter$}
    \end{axis}
\end{tikzpicture}
    	\caption{Projection onto the flow time and $N_f$ plane.}
    	\label{fig:RegressionTubeFVFlowProjection}
    \end{subfigure}
    \caption{Multi-variate linear regression of $\timing^f$ in function of $\totalProblemIter^c$ and $\totalProblemIter^f$ for the tube case with the FV-FE framework.}
    \label{fig:RegressionFlow}
\end{figure}
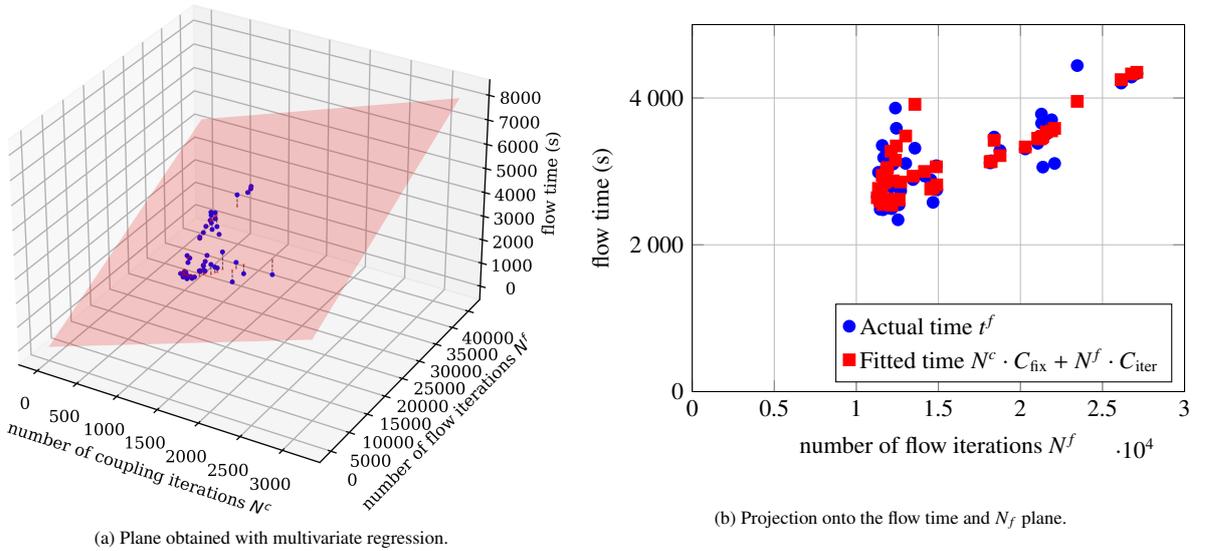

For the flow and solid time a multivariate linear regression is applied with two independent variables: the number of coupling iterations $\totalCoupleIter$ and the number of {\subproblemIter s} $\totalProblemIter^f$ or $\totalProblemIter^s$.
Neglecting any cost that occurs only once per simulation,
a zero intercept is considered, such that the result of the regression is effectively a plane through the origin.
As an example,
\Fig{RegressionTubeFVFlow3D} shows this plane for the flow time of the tube case solved with the FV-FE framework, 
while
\Fig{RegressionTubeFVFlowProjection} projects the same data into a plot over the flow solver iterations. 

The accuracy of the fit is assessed by the relative root mean square error (RRMSE), which is calculated for the flow time as
\begin{equation}
    \mathrm{RRMSE} = \sqrt{\frac{\sum^\numCalc\abs{
    \timing^f-(\totalProblemIter^c \cdot \costFix^f + \totalProblemIter^f \cdot \costIter^f)
    }^2}{\sum^m\abs{\timing^f}^2}},
    \label{equ:RRMSE}
\end{equation}
where $\numCalc$ is the number of calculations considered.
For the flexible tube with FV-FE the RRMSE of the flow time regression is 7.43 \%. This value shows that the regression is reasonably accurate.
The remaining difference is caused by natural variation in run time due to the varying loading and clock speed of the processors
as well as minor effects not considered in the model, like varying cost of solving the assembled matrix-vector system.
Also, optimization actions taken by a solver behind the scenes such as load balancing or optimization of certain solver parameters, e.g., relaxation factors, contribute to the variations in run time. 
This is especially the case for the commercial FV solver in the FV-FE framework\footnote{Accordingly, the RRMSE of the FE-FE framework, which uses only scientific code, is much lower for both test cases, see \Tab{RRMSE_Measures}.}.
For completeness, the RRMSE of the solid time regression is 3.44 \%.

\begin{figure}
	\centering
	\begin{tikzpicture}
    \smallIfElsevier 
    \begin{axis}[
        width=0.5\textwidth,
        height=0.4\textwidth,
        grid=both,
        xmin=0,
        xmax=2500,
        ymin=0,
        ymax=200,
        xlabel=number of coupling iterations $\totalProblemIter^c$,
        ylabel=coupling time ($\qty{}{\second}$),
        /pgf/number format/.cd,
        fixed,
        1000 sep={\ },
        legend style={
            at={(0.975,0.025)},
            anchor=south east,
            legend cell align=left,
            },
        ]
        \addplot [only marks,thick,blue,mark=*] table[x=it,y=time,col sep=comma]{Data/regression_coupling.csv};
        \addlegendentry{Actual time $\timing^c$}
        \addplot [only marks,thick,red,mark=square*] table[x=it,y=eq_time,col sep=comma]{Data/regression_coupling.csv};
        \addlegendentry{Fitted time $\totalProblemIter^c \cdot \costSolver^c$}
        \addplot[
		domain = 0:2500,
		samples = 10,thick,brown!60!black,dashed,] 
		{0.07951458020914454*x};
    \end{axis}
\end{tikzpicture}
	\caption{Linear regression of $\timing^c$ in function of $\totalProblemIter^c$ for the tube case with the FV-FE framework.}
	\label{fig:RegressionTubeFVCoupling}
\end{figure}
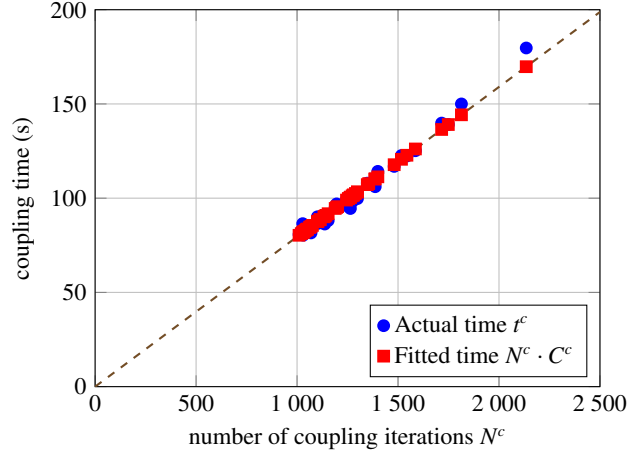

For the coupling time a linear regression analysis is performed with only one independent variable: the total number of coupling steps $\totalCoupleIter$.
Again, a zero intercept is assumed, such that the result is a line through the origin.
\Fig{RegressionTubeFVCoupling} shows the line for the tube case solved with the FV-FE framework.
The RRMSE value is 2.65 \%, which is again low, indicating that the assumed 
proportionality between coupling iterations and coupling cost is accurate.

As discussed in \Sec{convergence}, all cost coefficients scaling with the number of coupling iterations can be combined into $\costCoupleSum$.
Plugging these coefficients into \Equ{equivalentTime} gives a cost measure for the run time, which will be called the \textit{equivalent time}.
For the flexible tube case with FV-FE, the absolute value of the relative difference between the actual run time and the equivalent time is on average 3.80 \% and never higher than 14.62 \%. The first measure is referred to as the mean absolute percentage error (MAPE) and, analogously, the second is called the maximum absolute percentage error (maxAPE). These are defined as
\begin{equation}
    \mathrm{MAPE} = \frac{1}{\numCalc}\sum^\numCalc\abs{\frac{\timeSimulation-\costSimulation}{\timeSimulation}}
    \quad\mathrm{and}\quad
    \mathrm{maxAPE} = \max_\numCalc\left(\abs{\frac{\timeSimulation-\costSimulation}{\timeSimulation}}\right).
    \label{equ:MAPEmaxAPE}
\end{equation}

\begin{figure}
	\centering
	\begin{tikzpicture}
    \smallIfElsevier 
    \begin{axis}[
        width=0.5\textwidth,
        height=0.4\textwidth,
        grid=both,
        xmin=0,
        xmax=2500,
        ymin=0,
        ymax=8000,
        xlabel=number of coupling iterations $\totalCoupleIter$,
        ylabel=time ($\qty{}{\second}$),
        /pgf/number format/.cd,
        fixed,
        1000 sep={\ },
        legend style={
            at={(0.975,0.025)},
            anchor=south east,
            legend cell align=left,
            },
        ]
        \addplot [only marks,thick,blue,mark=*] table[x=c_it,y=time,col sep=comma]{Data/comparison_with_old_measure.csv};
        \addlegendentry{Actual time $\timeSimulation$}
        \addplot [only marks,thick,red,mark=square*] table[x=c_it,y=eq_time,col sep=comma]{Data/comparison_with_old_measure.csv};
        \addlegendentry{Equivalent time $\costSimulation$}
        \addplot [only marks,thick,green!60!black,mark=triangle*] table[x=c_it,y=old_eq_time,col sep=comma]{Data/comparison_with_old_measure.csv};
        \addlegendentry{Measure in literature}
        \addplot[
		domain = 0:2500,
		samples = 10,thick,brown!60!black,dashed,] 
		{3.098512158640545*x};
    \end{axis}
\end{tikzpicture}
	\caption{Difference between the actual time, the equivalent time used here, and the cost measure considering a constant cost per coupling iteration for the flexible tube case with the FV-FE framework. Note that the equivalent time with constant cost per coupling iteration forms a line through the origin.}
	\label{fig:EquivalentTime}
\end{figure}
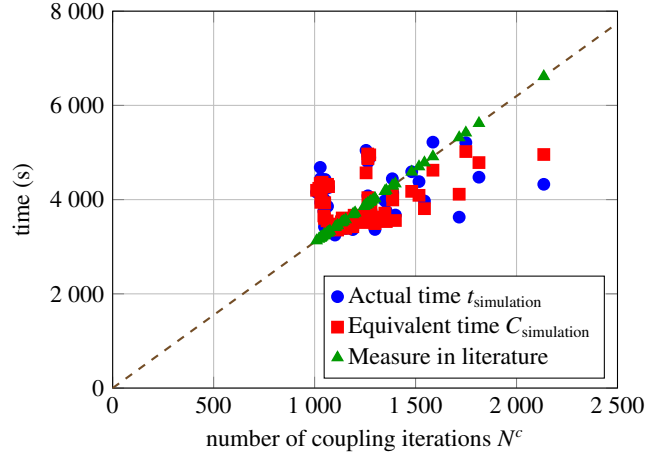

For comparison, assuming that the total cost scales only with the coupling iterations as common in literature,
i.e., setting $\costFac{iteration}^{f,s}=0$,
the deviation increases to an average of 13.58 \% and a maximum of 52.96 \%.
\Fig{EquivalentTime} illustrates the difference between these measures and their error with respect to the actual cost of the simulation.
Hence, the measure proposed in this work is not perfect, but definitely 
an improvement over the one predominant in literature.

An additional benefit of this equivalent time is that it is not affected by the occurring variations in run time, as the random variations present in the actual run time are filtered out.
For this reason, the equivalent time will be used as cost measure instead of the actual time from here on.
\Tab{RegressionSummary} gives a summary of the cost factors and the differences between the actual and equivalent run time for the lid-driven cavity and flexible tube case with both frameworks
\footnote{
At first, it might seem strange that, for the lid-driven cavity case in \Tab{RegressionSummary}, a fixed-point iteration in the FV flow solver is
more expensive than a Newton step in the flow solver of its FE-FE counterpart.
However,
the two frameworks did not only use different software, but were also run on different hardware infrastructure,
prohibiting a direct comparison.}
%

\begin{table}
	\caption{Cost factors in seconds per corresponding iteration and difference between the actual and equivalent run time, expressed with the mean and maximum absolute percentage error.}
	\label{tab:RegressionSummary}
	\begin{center}
        \smallIfElsevier 
		\begin{tabular}{l|*{5}{l}|l|ll}
			& $\costFix^f$ & $\costIter^f$ & $\costFix^s$ & $\costIter^s$ & $\costCouple$ & $\costCoupleSum$ & MAPE & maxAPE \\
			\hline
			\textbf{Lid-driven cavity} & & & & & & & & \\
			$\quad$ FE-FE     & 0.0208 & 0.0298 & 0.0002 & 0.0023 & 0.0004 & 0.0214 & 0.33 \% & 0.89 \% \\
			$\quad$ FV-FE     & 0.8985 & 0.0727 & 0.0243 & 0.0041 & 0.0218 & 0.9446 & 3.86 \% & 15.76 \% \\
            \textbf{Flexible tube} & & & & & & & & \\
			$\quad$ FE-FE     & 0.6459 & 1.4756 & 0.0128 & 0.2076 & 0.1873 & 0.8460 & 0.22 \% & 0.58 \% \\
			$\quad$ FV-FE     & 1.1542 & 0.1068 & 0.1587 & 0.2510 & 0.0795 & 1.3924 & 3.80 \% & 14.62 \% 
		\end{tabular}
	\end{center}
\end{table}

The methodology explained in this section was illustrated by the flexible tube simulations with the FV-FE framework. Quality measures of the regression for other combinations of problems and frameworks are given in \App{fullRegression}. 


\subsection{Results of parameter study}
\label{sec:parameterStudy}

The primary research question of this work is how the overall performance of the partitioned FSI algorithm is influenced by limiting the number of \subproblemIter s per solver call.
%
Towards this goal,
simulations were run multiple times with varying values for $\iterPerCall^f$ and $\iterPerCall^s$. 
By virtue of the new convergence criterion, each of these runs solves the coupled problem up to the same tolerance.
Two test cases investigated by two FSI software frameworks led to a total of four parameter studies.
The following subsections present their results, focusing
on how varying $\iterPerCall^f$ and $\iterPerCall^s$ influences different quantities and measures.

\subsubsection{Iteration counts} \label{sec:Results_IterCounts}

\begin{figure}[ht]
    \centering
    \begin{subfigure}{.33\textwidth}
    	\includegraphics[width=\textwidth,trim={0 55 0 5},clip]{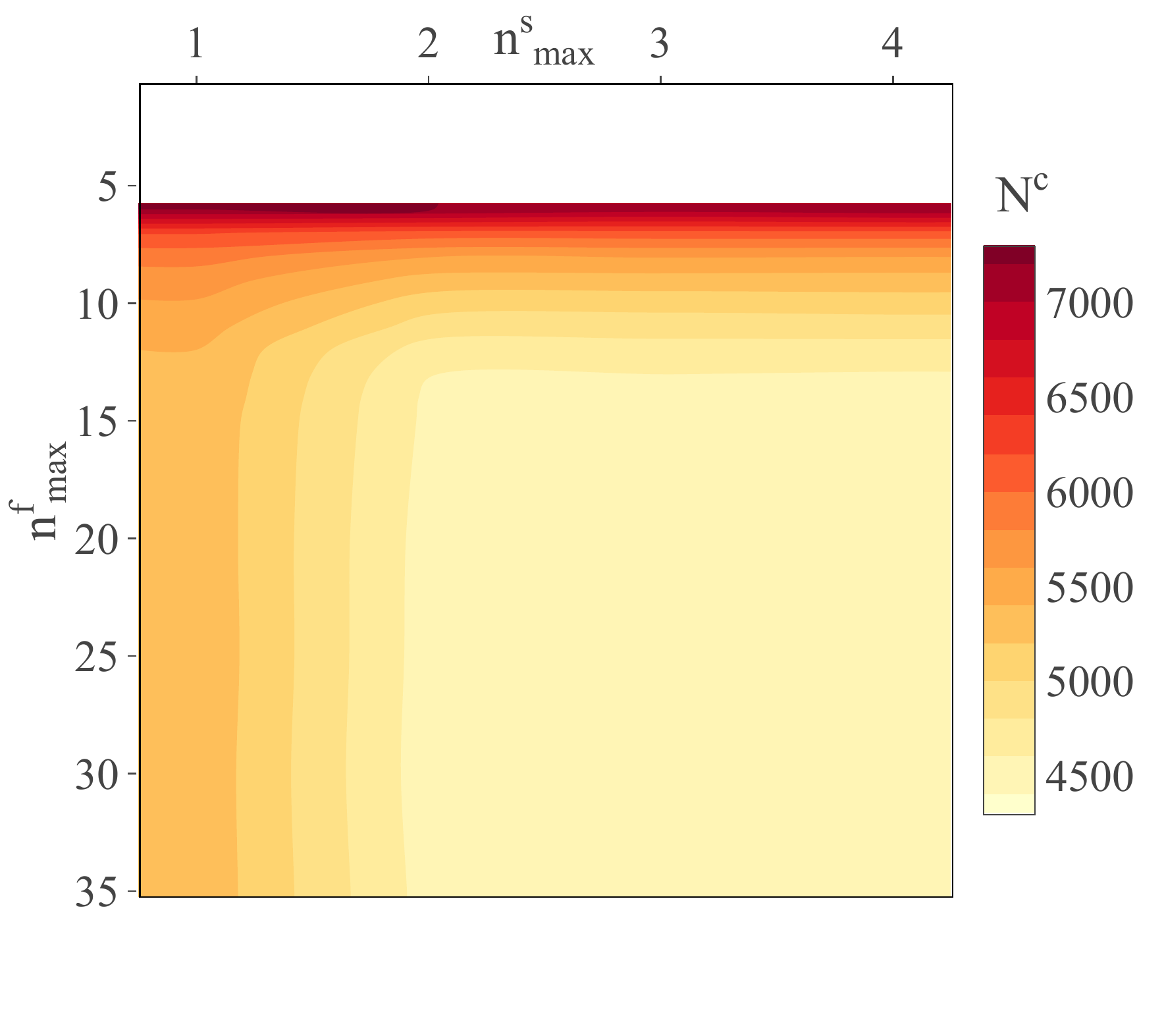}
        \caption{Coupling iterations $\totalCoupleIter$.}
    \end{subfigure}
    \begin{subfigure}{.33\textwidth}
        \includegraphics[width=\textwidth,trim={0 55 -14 5},clip]{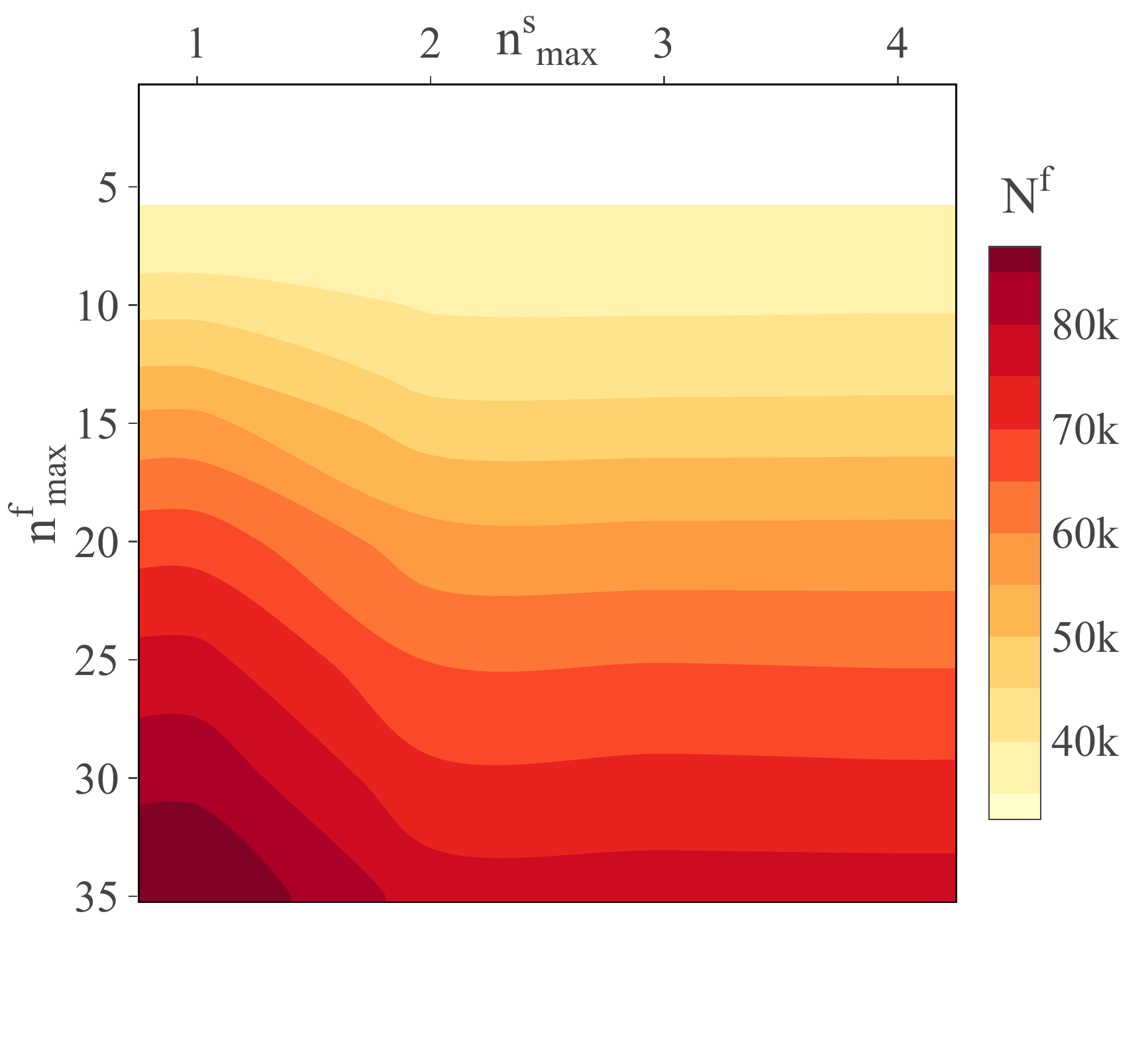}
    	\caption{Flow \subproblemIter s $\totalProblemIter^f$.}
    \end{subfigure}
    \begin{subfigure}{.33\textwidth}
        \includegraphics[width=\textwidth,trim={0 55 -14 5},clip]{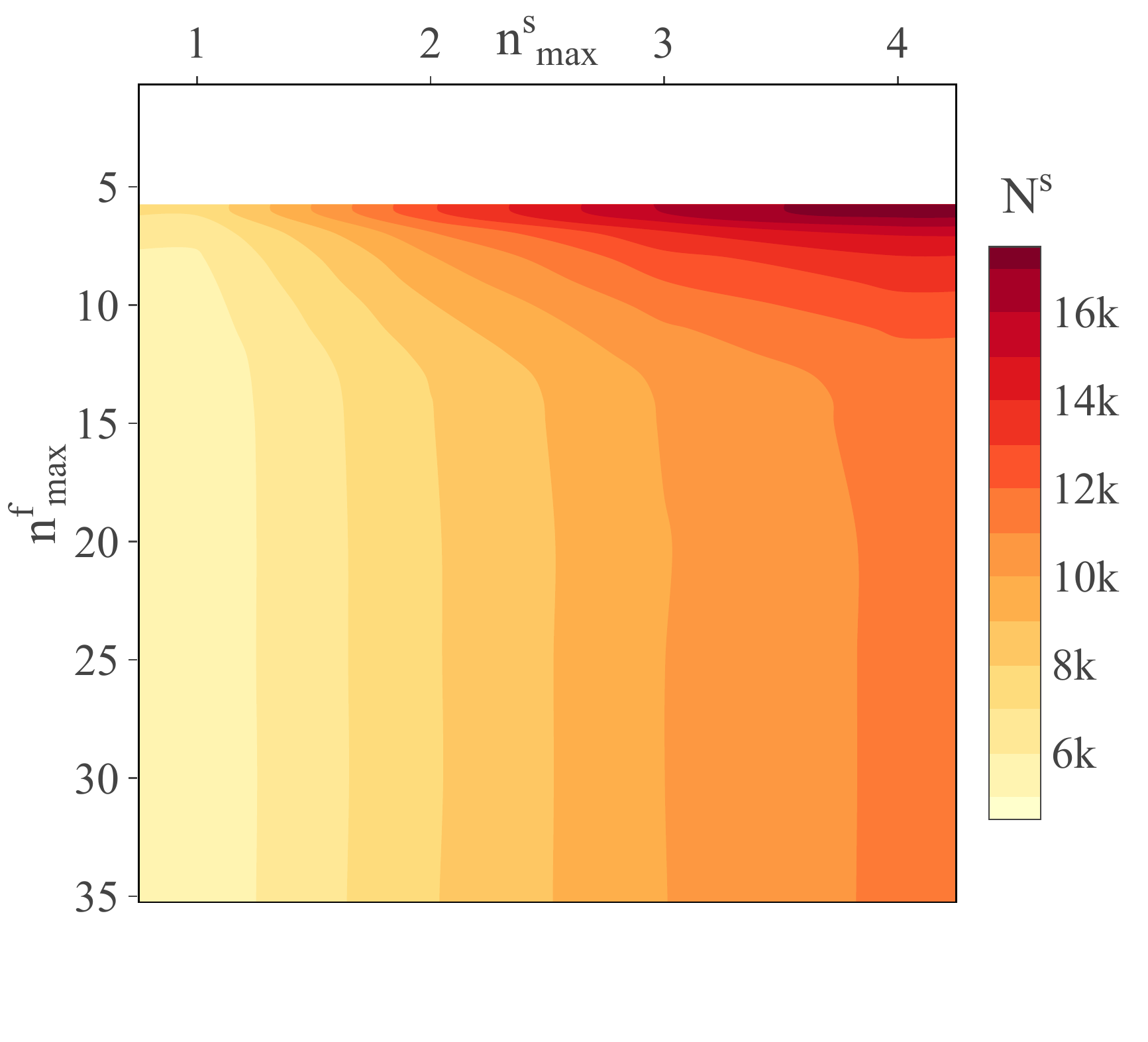}
    	\caption{Solid \subproblemIter s $\totalProblemIter^s$.}
    \end{subfigure}
    \caption{Different iteration counts plotted over $\iterPerCall^f$ and $\iterPerCall^s$ for the lid-driven cavity case simulated with the FV-FE framework.}
    \label{fig:ContourPlots_FV_Cavity}
\end{figure}

\begin{figure}[ht]
    \centering
    \begin{subfigure}{.33\textwidth}
    	\includegraphics[width=\textwidth,trim={0 55 -14 5},clip]{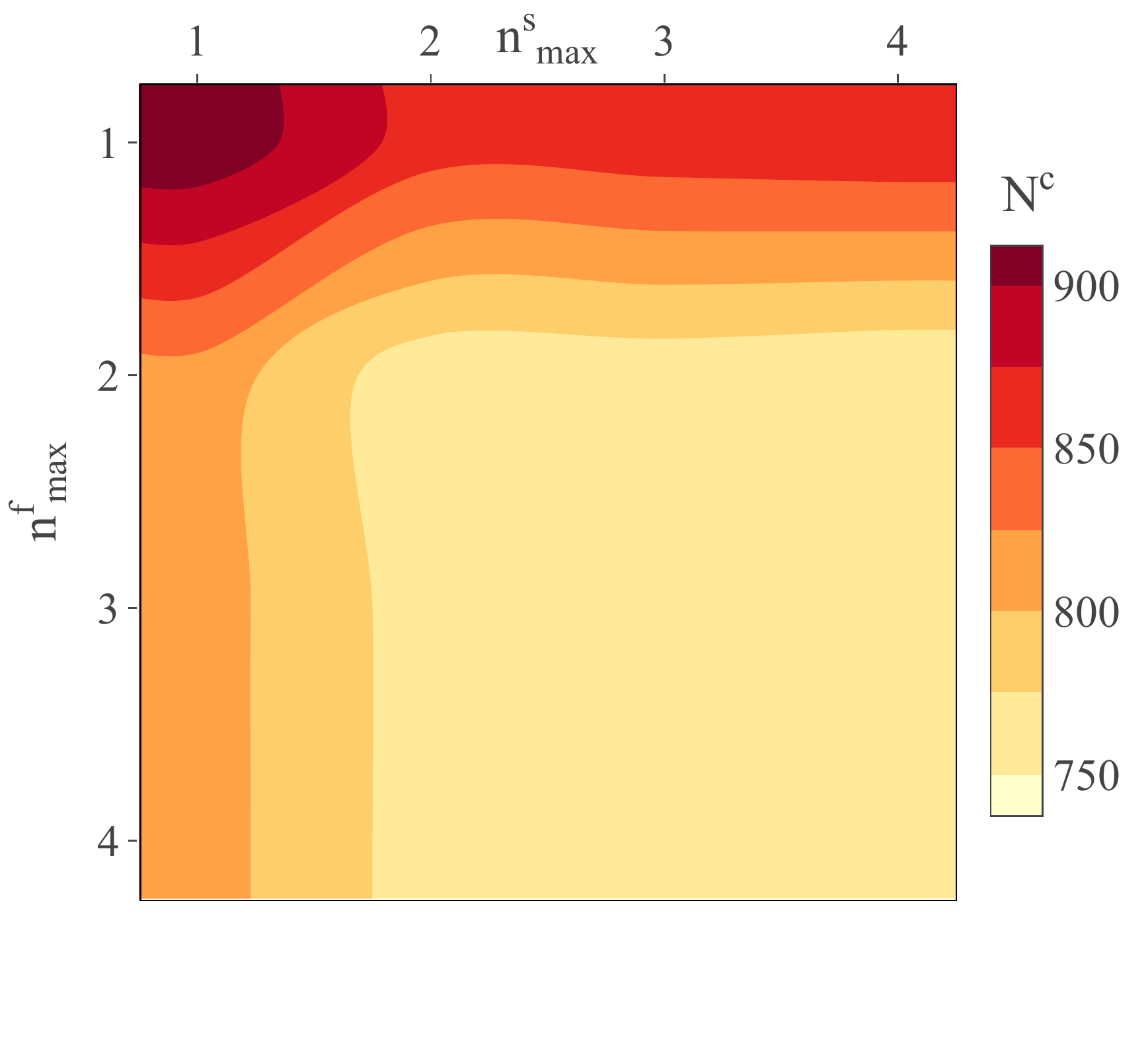} 
        \caption{Coupling iterations $\totalCoupleIter$.}
    \end{subfigure}
    \begin{subfigure}{.33\textwidth}
        \includegraphics[width=\textwidth,trim={0 55 0 5},clip]{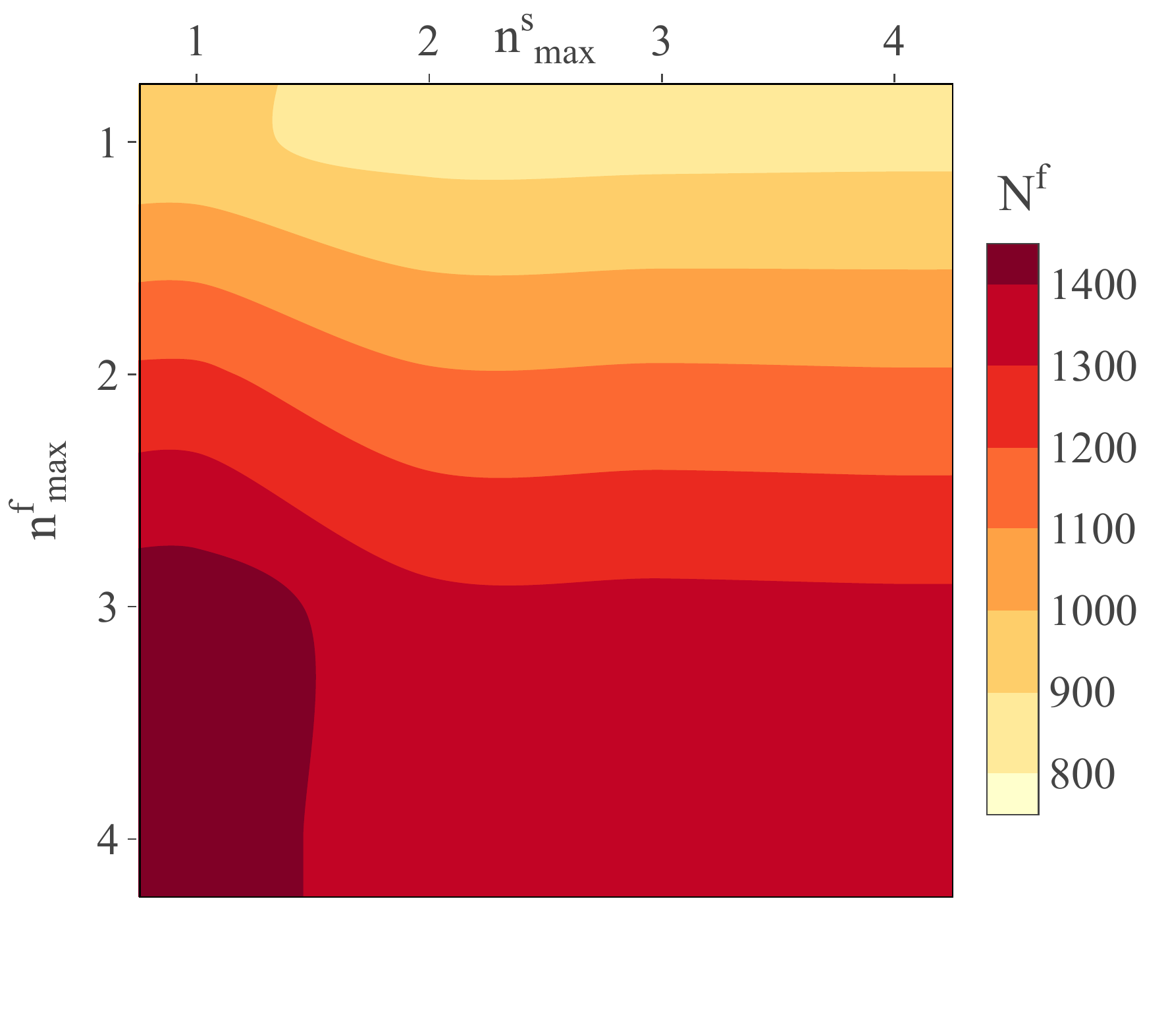}
    	\caption{Flow \subproblemIter s $\totalProblemIter^f$.}
    \end{subfigure}
    \begin{subfigure}{.33\textwidth}
        \includegraphics[width=\textwidth,trim={0 55 0 5},clip]{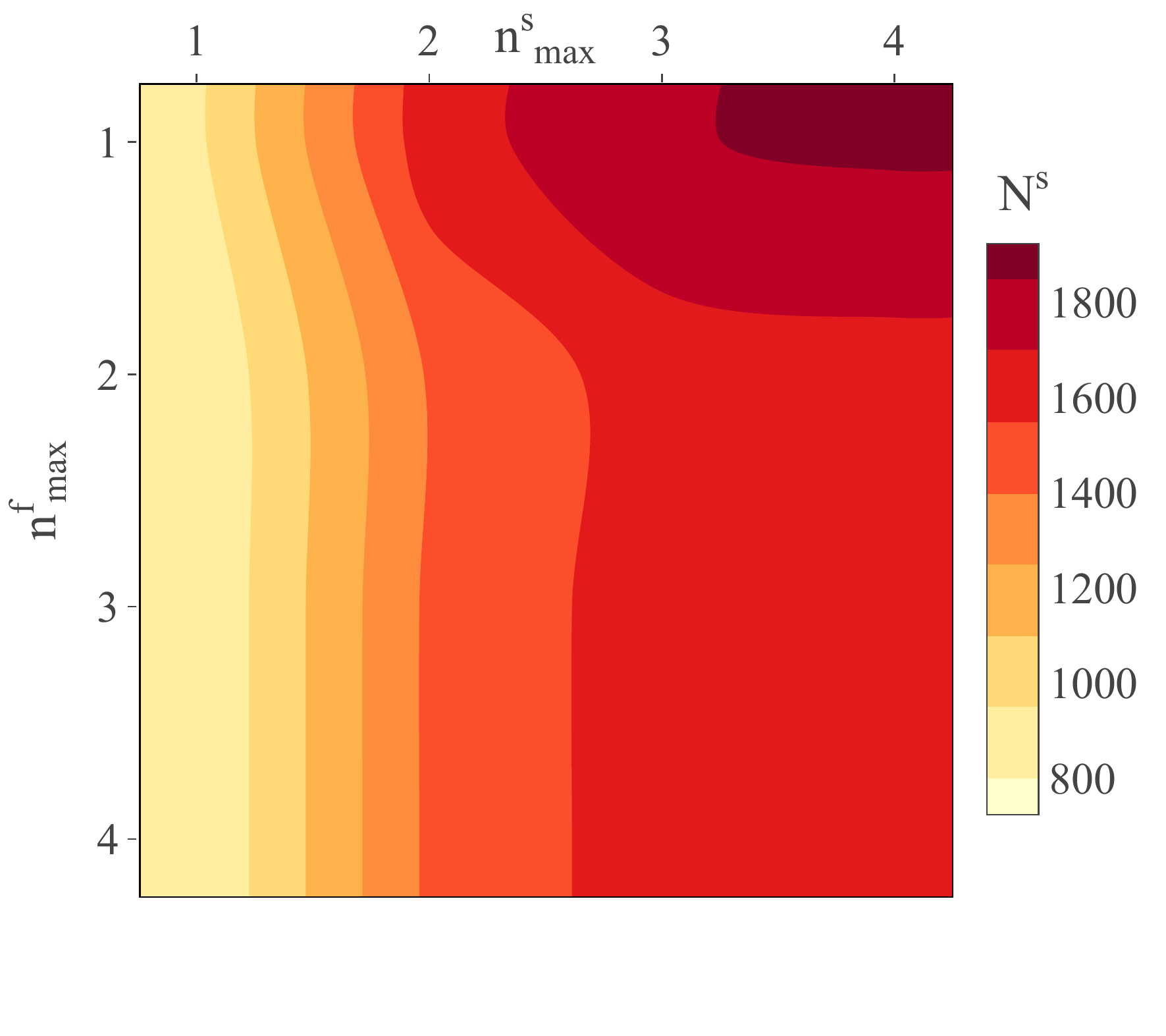}
    	\caption{Solid \subproblemIter s $\totalProblemIter^s$.}
    \end{subfigure}
    \caption{Different iteration counts plotted over $\iterPerCall^f$ and $\iterPerCall^s$ for the tube case simulated with the FE-FE framework.}
    \label{fig:ContourPlots_FE_Tube}
\end{figure}

One central aspect is the impact of the \subproblemIter s per solver call on the three global iteration counts, i.e., 
the total number of coupling iterations $\totalCoupleIter$ as well as
the total number of \subproblemIter s of the flow and solid solver, $\totalProblemIter^f$ and $\totalProblemIter^s$, respectively.
Therefore, \Fig{ContourPlots_FV_Cavity} and \Fig{ContourPlots_FE_Tube} illustrate their dependence on $\iterPerCall^f$ and $\iterPerCall^s$ in three contour plots each.
To increase the representativity of the examples,
\Fig{ContourPlots_FV_Cavity} is based on the lid-driven cavity case simulated with the FV-FE framework,
while \Fig{ContourPlots_FE_Tube} visualizes the tube case for the FE-FE framework.

The overall trends observed for these two examples are in very good agreement. The biggest difference is that the
plots obtained for the FV flow solver are more 'fine-grained',
since it relies on fixed-point rather than Newton steps, requiring more but typically less expensive iterations.
The key findings of the three subplots are:
\begin{enumerate}[(a)]
\item The number of \textit{coupling iterations} clearly decreases when running more \subproblemIter s per solver call, so that iterating to full convergence in every call yields the smallest $\totalCoupleIter$.
This observation, in line with the reasoning of \Sec{IterationsPerCall}, is expected
since the interface data fed back by the subproblem into the coupling loop is ensured to be as accurate as possible given the input data.
The other way around, reducing $\iterPerCall^f$ and/or $\iterPerCall^s$ deteriorates the quality of the exchanged data fields, leading to more coupling iterations.
This effect can even be severe enough to cause divergence of the coupling, as observed for the FV-FE framework in \Fig{ContourPlots_FV_Cavity} if $\iterPerCall^f<6$.
\item The main influence on the \textit{flow iterations} is that the total number $\totalProblemIter^f$ grows with $\iterPerCall^f$. 
This indicates that for increasing $\iterPerCall^f$, the increase in $\totalProblemIter^f$ resulting from the additional \subproblemIter s performed within each coupling step outweighs those saved due to the reduced number of coupling iterations, discussed in (a).
Keeping $\iterPerCall^f$ fixed reveals a secondary trend:
running more solid iterations per solver call decreases
the total number of flow iterations. Again, this is explicable by the improved quality of the interface data that is passed back by the solid solver and serves as the fluid's boundary condition.
\item The contour plot of the \textit{solid iterations} and its characteristics are principally the transpose of those discussed for the flow iterations in (b):
$\totalProblemIter^s$ primarily depends on $\iterPerCall^s$, but is influenced by $\iterPerCall^f$ as well.
\end{enumerate}

The contour plots of the two missing combinations, i.e., the lid-driven cavity with the FE-FE setup and  the tube case simulated with the FV-FE solver framework, 
are omitted for the sake of conciseness.
They show the same trends,
which suggest the discussion above is to some extent general, in that the underlying effects are not very problem-dependent, or at least observable for typical FSI simulations.

\subsubsection{Computational cost}

Investigating the iteration counts $\totalCoupleIter$, $\totalProblemIter^f$, $\totalProblemIter^s$ is certainly of scientific interest,
but in practice
the most important quantity to consider is
the simulation's total run time,
in this work represented by the equivalent time measure 
defined in \Equ{equivalentTime}.
Since it is computed as a weighted sum of $\totalCoupleIter$, $\totalProblemIter^f$, and $\totalProblemIter^s$,
all effects discussed in the previous section affect
the equivalent time measure too.
Their significance, however, is determined by the 
weighting factors, i.e., the cost of one coupling, flow, or solid iteration approximated by the regression model.
Unlike the more general dependencies illustrated in \Fig{ContourPlots_FV_Cavity} and \Fig{ContourPlots_FE_Tube},
the equivalent time measure therefore strongly depends on other aspects such as the problem itself, the solver framework, the HPC architecture, and so on.

\begin{figure}[ht]
    \centering
    \begin{subfigure}{.49\textwidth}
    	\includegraphics[width=\textwidth,trim={0 55 0 0}, clip]{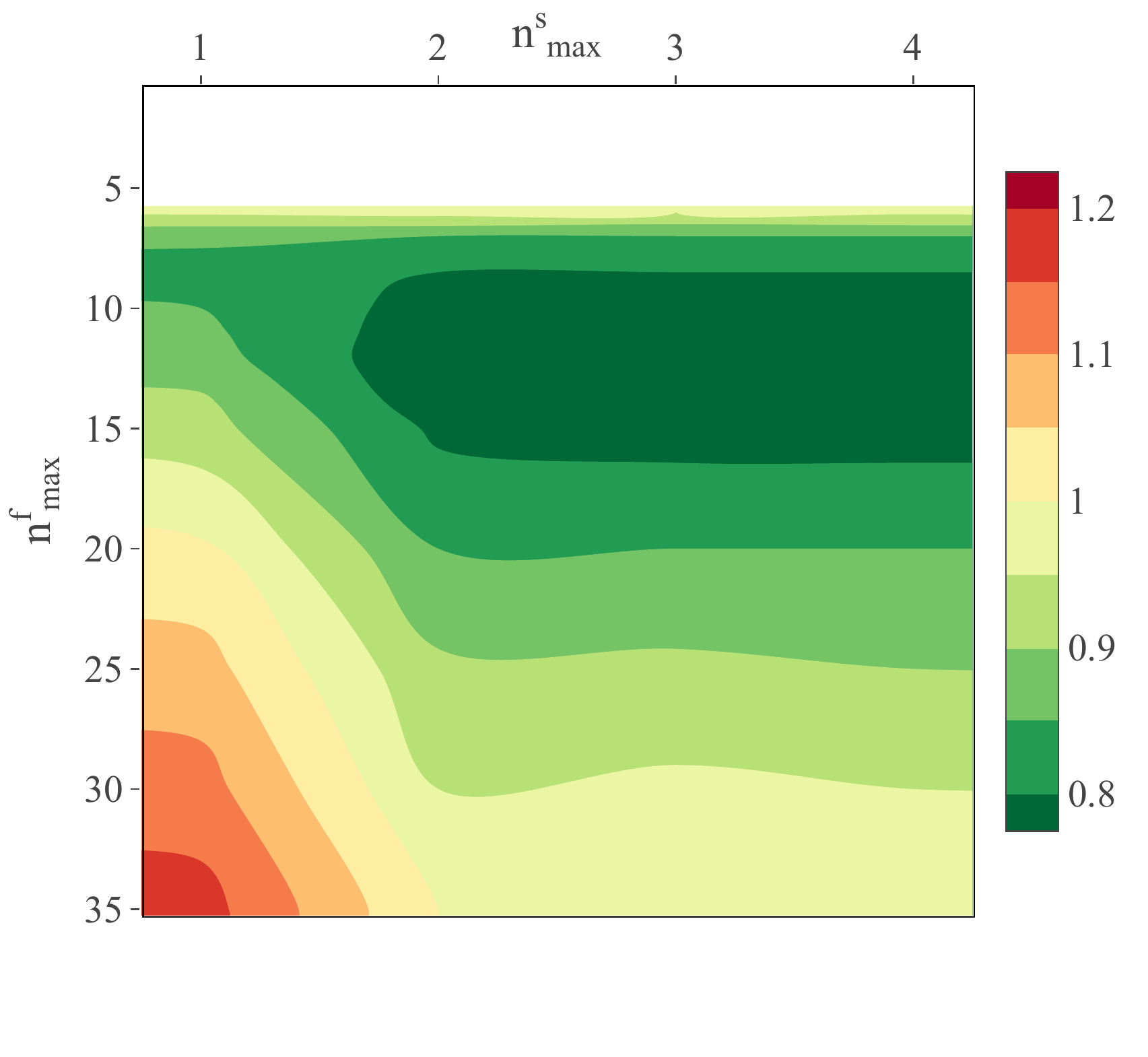} 
        \caption{FV-FE framework.}
    \end{subfigure}
    \begin{subfigure}{.49\textwidth}
        \includegraphics[width=\textwidth,trim={0 44 0 0},clip]{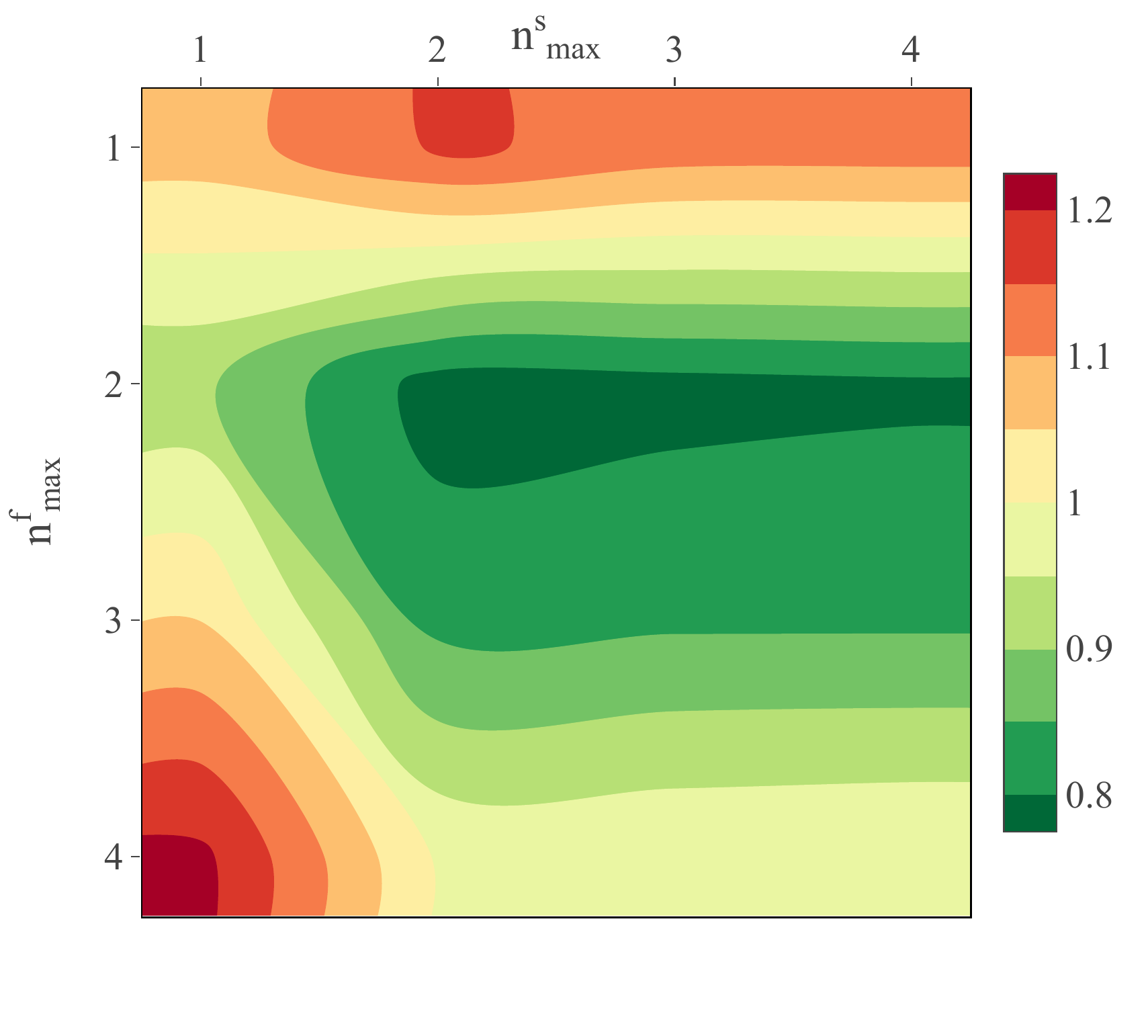}
    	\caption{FE-FE framework.}
    \end{subfigure}
	\caption{Contour plots of the equivalent time measure for the lid-driven cavity case, normalized with respect to the value obtained for full convergence.} \label{fig:ContourPlots_EqvCavity_Time}
\end{figure}

For the lid-driven cavity case, the parameter study of the two solver frameworks results in the contour plots
depicted in \Fig{ContourPlots_EqvCavity_Time}. To make them more comparable, the
values are normalized by dividing through the equivalent time obtained for iterating each solver call to full convergence.
The most striking observation is that,
despite resulting from different solver frameworks run on different computer systems,
the overall characteristics of the plots are very similar.
The only major difference is found for the top part, i.e., for small values of $\iterPerCall^f$.
While the FV-FE framework diverged for $\iterPerCall^f<6$ as discussed in the previous section,
no lower bound is apparent for the FE-FE setup,
as performing a single Newton iteration per flow solver call is already sufficient to ensure convergence. Nevertheless, the negative impact of $\iterPerCall^f=1$ on the stability leads to
more coupling iterations and therefore higher computational cost.

Apart from this,
the trends apparent in the two plots are in good agreement.
They clearly show that the total run time is not minimized by iterating to full convergence in every solver call.
Accordingly, the optima of the run time and the number of coupling iterations $\totalCoupleIter$ do not match.
Moreover, the least efficient choice in both cases is to set $\iterPerCall^s=1$
and $\iterPerCall^f=\infty$.
%
The area of lowest computational cost,
on the other hand,
is found for some mid-range values for both parameters, like $\iterPerCall^f=\iterPerCall^s=2$ for FE-FE
or $\iterPerCall^f \approx 12$ and $\iterPerCall^s=2$ for the FV-FE framework.
While further increasing $\iterPerCall^f$ reduces the efficiency significantly, adding more solid iterations in this case does not have a big impact, since solving the solid problem is much cheaper for both setups.
%

Although the most significant speed up with respect to full convergence in every solver call, i.e., $\iterPerCall^f=\iterPerCall^s=\infty$,
is only about $22 \%$,
the parameter study of the lid-driven cavity case clearly shows that the number of \subproblemIter s per solver call
has a significant influence on the computational cost.
In order to be efficient, an FSI solver framework should therefore take
this impact into account. \\

\begin{figure}[ht]
    \centering
    \begin{subfigure}{.49\textwidth}
    	\includegraphics[width=\textwidth,trim={0 55 0 0}, clip]{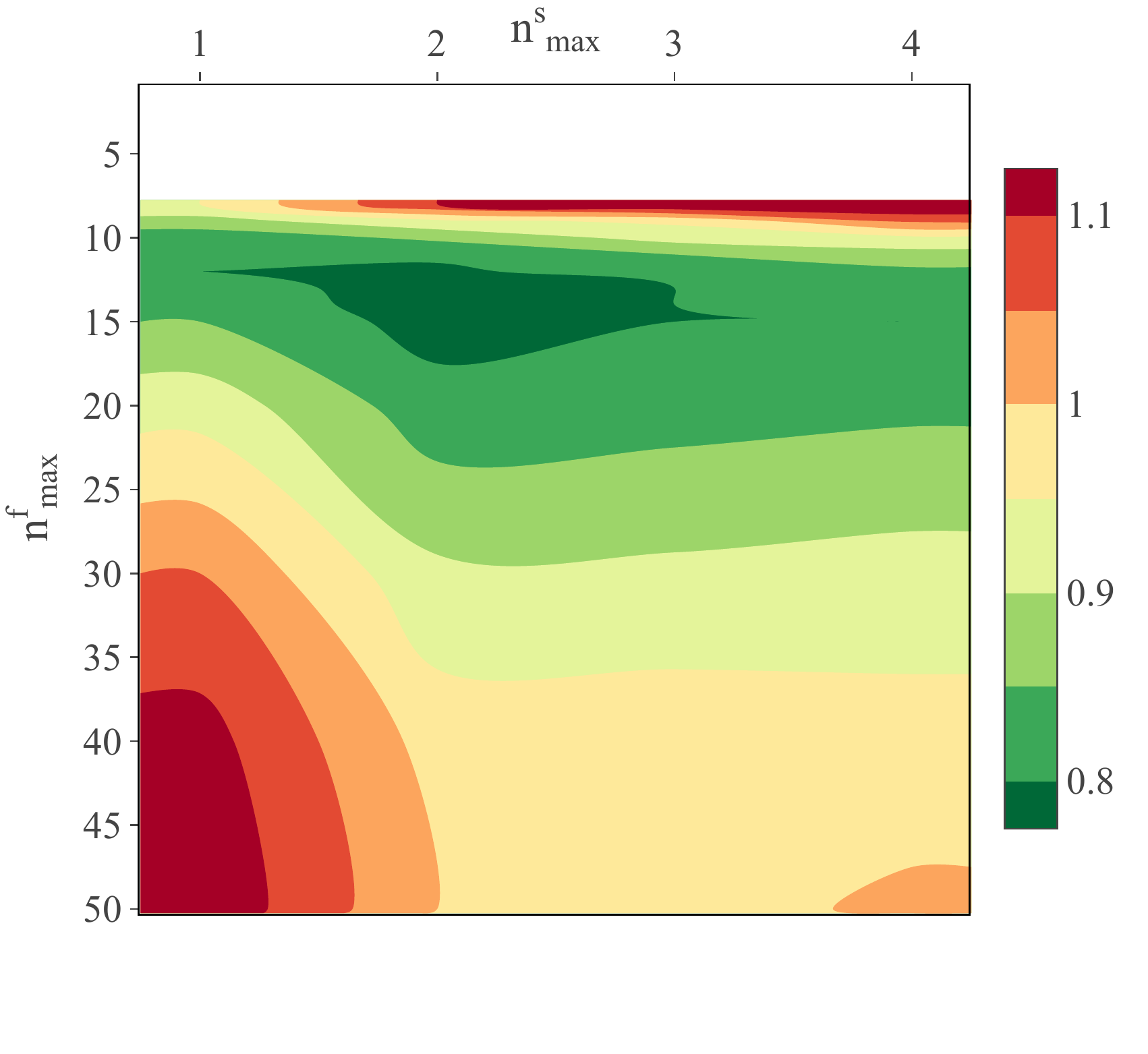} 
        \caption{FV-FE framework.}
    \end{subfigure}
    \begin{subfigure}{.49\textwidth}
        \includegraphics[width=\textwidth,trim={0 44 0 0},clip]{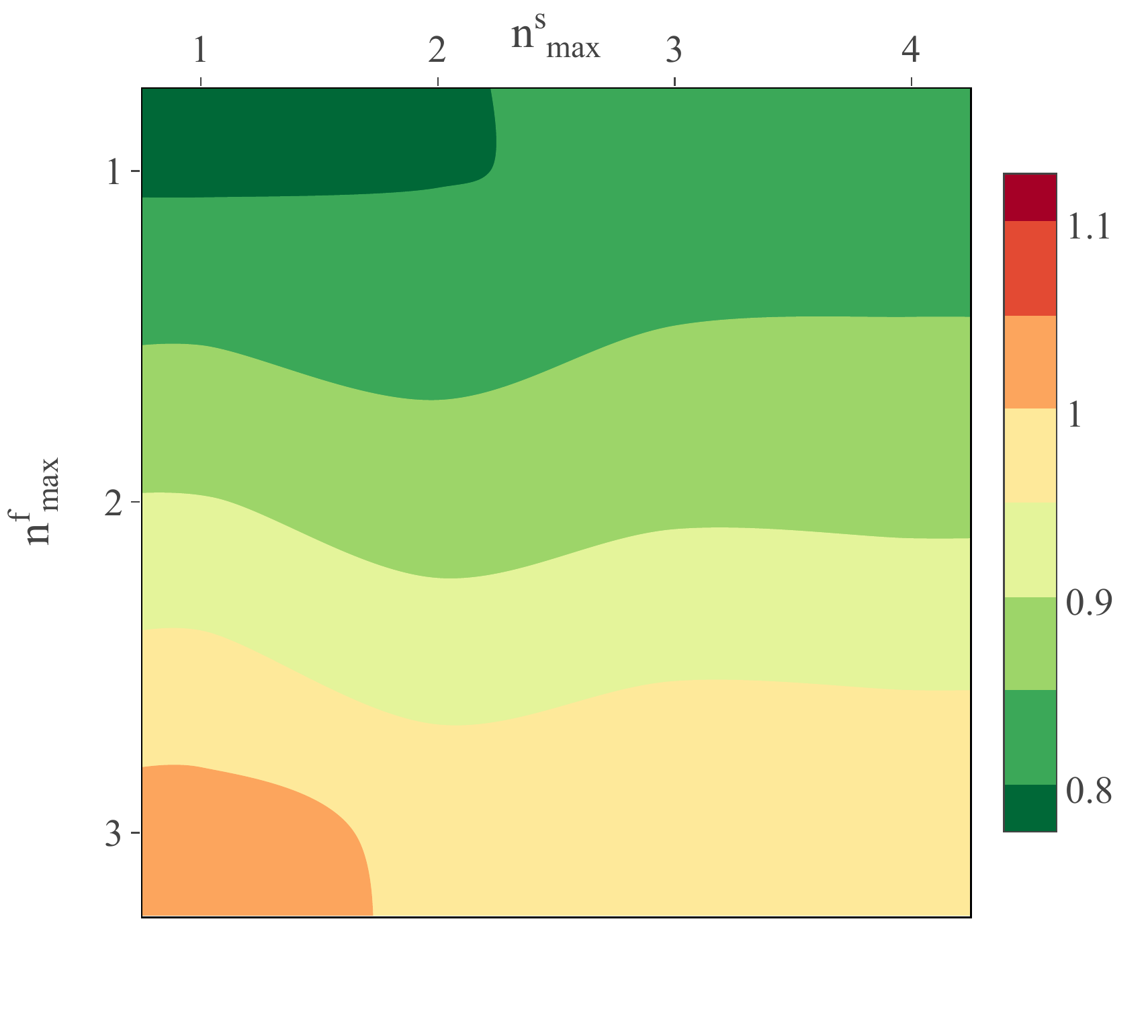}
    	\caption{FE-FE framework.}
    \end{subfigure}
	\caption{Contour plots of the equivalent time measure for the flexible tube case, normalized with respect to the value obtained for full convergence.}
    \label{fig:ContourPlots_Tube_EqvTime}
\end{figure}

Analogously to the discussion of the lid-driven cavity, \Fig{ContourPlots_Tube_EqvTime} illustrates the computational cost as a function of $\iterPerCall^f$ and $\iterPerCall^s$ for the flexible tube case.
Undoubtedly, the difference between the two solver frameworks is more pronounced than in the previous figure.
Nonetheless, the biggest difference is again observed for very small values of $\iterPerCall^f$.
Like in the previous case, the FV flow solver requires a certain minimum number of \subproblemIter s, here $\iterPerCall^f\gtrsim8$, to prevent divergence.
For the FE-FE framework, on the other hand, one Newton iteration per call is not only enough to converge,
but in this case even results in the lowest computational cost.
A likely explanation is that one Newton iteration is already accurate enough to yield similar results as in the green region observed in the FV-FE framework around 
$\iterPerCall^f \approx 12$, so that the part above with slightly slower convergence is
never reached for the FE-FE setup.
This reasoning is supported by the observation that for all FE flow solver calls of the parameter study, a maximum of three Newton iterations was enough to converge, 
indicating that Newton steps are a very effective method of handling the flow problem's nonlinearity for this test case.
%

In both plots,
increasing the number of flow iterations per solver call leads to a smooth yet significant growth in computational cost, while the solid iterations have a far lower impact.
Furthermore, setting $\iterPerCall^s=1$ and $\iterPerCall^f=\infty$ is the most expensive choice for both frameworks.
In that regard, the two plots are similar. 
However, this maximum is much more distinct in the FV-FE case than for FE-FE.
Given the three contour plots in \Fig{ContourPlots_FE_Tube},
it is clear that the peak in the lower left corner stems from 
the flow solver's cost.
As the cost factors listed in \Tab{RegressionSummary}
do not reveal the flow solver to be more expensive in the FV-FE framework
relative to the other costs than for FE-FE,
the difference is likely to be caused by a locally higher 
increase of the number of flow iterations $\totalProblemIter^f$.
This hypothesis is supported by the tables presented in \App{fullResults}.

\subsubsection{Further remarks}

The contour plots discussed so far illustrated the most important trends and findings of the parameter studies.
Further, \App{fullResults} contains four tables,
%
each corresponding to a particular parameter study, i.e., one test case and framework.
They list the coupling, flow, and solid iterations as well
as the relative equivalent times obtained for all runs belonging to that particular parameter study.
%
The following additional observations and remarks are noteworthy:
\begin{itemize}
    \item For both test cases, the FE-FE framework requires slightly fewer coupling iterations. 
    This effect is caused by the different types of \subproblemIter s.
    As the residual reduction is much more significant for a Newton step than it is for a fixed-point iteration,
    the first subproblem residual satisfying the convergence tolerance 
    will typically be considerably lower for the FE solver than for its FV counterpart.
    Consequently, even for the same residual tolerance the FE results fed back into the coupling loop are on average more accurate.
    In line with the reasoning in \Sec{IterationsPerCall} and \Sec{Results_IterCounts}, this leads to a lower $\totalCoupleIter$.
    The effect is further amplified since each solver call includes at least one {\subproblemIter}, i.e., a Newton step for FE or a fixed-point iteration for FV, even if the subproblem has already converged.
    While the higher number of coupling iterations is not problematic in itself, it results in a higher number of modes in the quasi-Newton coupling technique which risk to be (almost) linear dependent. This requires an efficient filtering technique, especially when the reuse parameter is high.

    Numerical experiments showed, however, that the higher number of coupling iterations observed for the FV solver can be countered by performing batches of fixed-point \subproblemIter s, e.g., groups of ten, and checking convergence of the flow solver only after each batch.

    
    \item In literature,
    it is common to monitor the coupling's convergence via the fixed-point residual of the interface displacement $\Rk{k}$, defining either an absolute or relative tolerance.
    To give an idea of how strict the convergence criteria were chosen for the parameter study, \Tab{RkDx} lists for all four studies both the average and maximum value of $\normE{\Rk{k}}$ and $\frac{\normE{\Rk{k}}}{\normE{\vect{d}^k}}$, determined upon convergence for each time step. 
    %
    \item As explained in \Rem{Increment}, the convergence criterion is only exact if the IQN update increment vanishes. Therefore, \Tab{RkDx} provides the average and maximum norm of the IQN increment norm at the end of a time step for all studies,
    confirming this minor aspect can safely be ignored in practice.
    \item By virtue of the new convergence criterion, all runs within a parameter study are converged up to the same tolerance. In this way, their results are virtually the same.
%
    To quantify the deviation, after each time step of each run in the parameter study, the computed interface displacement field $\vect{d}^k$ is compared to
    that obtained for iterating to full convergence in every call $\vect{d}^k_{ref}$, by evaluating 
    \begin{align}
        \frac{\normE{\vect{d}^k-\vect{d}^k_{ref}}}{\sqrt{\nInterfaceDofs}} ~,
    \end{align}
    where $\nInterfaceDofs$ is the number of displacement degrees of freedom at the FSI interface.
    In its last two columns, \Tab{RkDx} list the average and maximum values of this deviation for each study.
\end{itemize}

\begin{table}
    \caption{Four different quantities that are evaluated after each time step of all simulations within the parameter study. The average and maximum value are reported. First, the norm of the fixed-point residual of the interface deformation $\Rk{k}$ is given, followed by its relative magnitude with respect to the interface deformation $\vect{d}^k$.
    The third quantity is the norm of the IQN increment $\Delta \vect{\tilde{d}}^k_{IQN}$,
    while the last one evaluates the deviation from the results obtained for $\iterPerCall^f=\iterPerCall^s=\infty$.} \label{tab:RkDx}
	\begin{center}
        \smallIfElsevier 
		\begin{tabular}{l | c c c c c c | c c}
			& \multicolumn{2}{c}{$\normE{\Rk{k}}$} & \multicolumn{2}{c}{$ \frac{\normE{\Rk{k}} }{ \normE{\vect{d}^k}}$} &\multicolumn{2}{c|}{$\normE{\Delta\vect{\tilde{d}}^k_{IQN}}$} & 
            \multicolumn{2}{c}{$\frac{\normE{\vect{d}^k-\vect{d}^k_{ref}}}{\sqrt{\nInterfaceDofs}}$}   \\
			& Average & Max & Average & Max & Average & Max & Average & Max \\
			\hline
			\textbf{Lid-driven cavity} & & & & & & & & \\
			$\quad$ FE-FE     & $1.98\E{-09}$ & $6.58\E{-07}$ & $7.92\E{-09}$ & $7.84\E{-06}$ & $8.53\E{-10}$ & $5.80\E{-07}$ & $5.12\E{-10}$ & $7.36\E{-9}$ \\
			$\quad$ FV-FE     & $2.54\E{-11}$ & $9.04\E{-10}$ & $1.18\E{-09}$ & $1.49\E{-06}$ & $4.95\E{-11}$ & $9.14\E{-10}$ & $3.99\E{-11}$ & $1.93\E{-10}$ \\
			\textbf{Flexible Tube} & & & & & & & & \\
			$\quad$ FE-FE     & $1.42\E{-11}$ & $1.90\E{-10}$ & $2.18\E{-08}$ & $1.10\E{-06}$ & $9.85\E{-12}$ & $1.56\E{-10}$ & $7.48\E{-12}$ & $2.63\E{-11}$ \\
			$\quad$ FV-FE     & $9.99\E{-12}$ & $2.33\E{-10}$ & $9.04\E{-09}$ & $5.44\E{-07}$ & $1.40\E{-11}$ & $2.56\E{-10}$ & $2.99\E{-12}$ & $8.34\E{-12}$ \\
		\end{tabular}
	\end{center}
\end{table}

The tables in \App{fullResults} show that the gain in computational efficiency and the location of the optimum depend on the specific problem and the selected solvers.
However, also the coupling technique and settings have a significant influence.
To illustrate this, 
\App{noReuse} presents
two additional parameter studies
for the flexible tube case with identical settings, but without reuse of data from past time steps in the IQN Jacobian approximation ($q$=0).
For both tables, fully converging the subproblems in every call leads to the worst performance.
The optimum is found by performing 1 Newton iteration in FE and around 8 fixed-point iterations in FV.
Note that for these choices the computational time is approximately halved with respect to the reference.

\subsubsection{Impact of cost factors}

The parameter studies demonstrate that limiting the number of \subproblemIter s per solver call has a significant influence on the computational cost of the partitioned simulation,
and that this cost is accurately represented by the new equivalent time measure, which relies on a weighted sum
of the
iteration counts $\totalCoupleIter$, $\totalProblemIter^f$, and $\totalProblemIter^s$.
%
%
Unfortunately, however, both weighting factors and iteration counts
are case-specific:
while the weighting factors refer to the cost of one iteration and hence in particular depend on the software framework and computer architecture,
the iteration counts are influenced for example by the chosen solver tolerances and convergence criteria.

Nevertheless, \Sec{Results_IterCounts} showed that these 
iteration counts follow rather general trends.
%
Taking the iteration counts obtained in the parameter studies
and artificially setting the cost factors $\costCoupleSum$, $\costIter^f$, and $\costIter^s$,
therefore allows to approximate the computational costs of different cost scenarios.

\begin{figure}[ht]
    \centering
    \captionsetup[subfigure]{,slc=off,margin={5pt,5pt}}
    \begin{subfigure}[t]{.33\textwidth}
        \includegraphics[width=\textwidth,trim={0 50 10 5},clip]{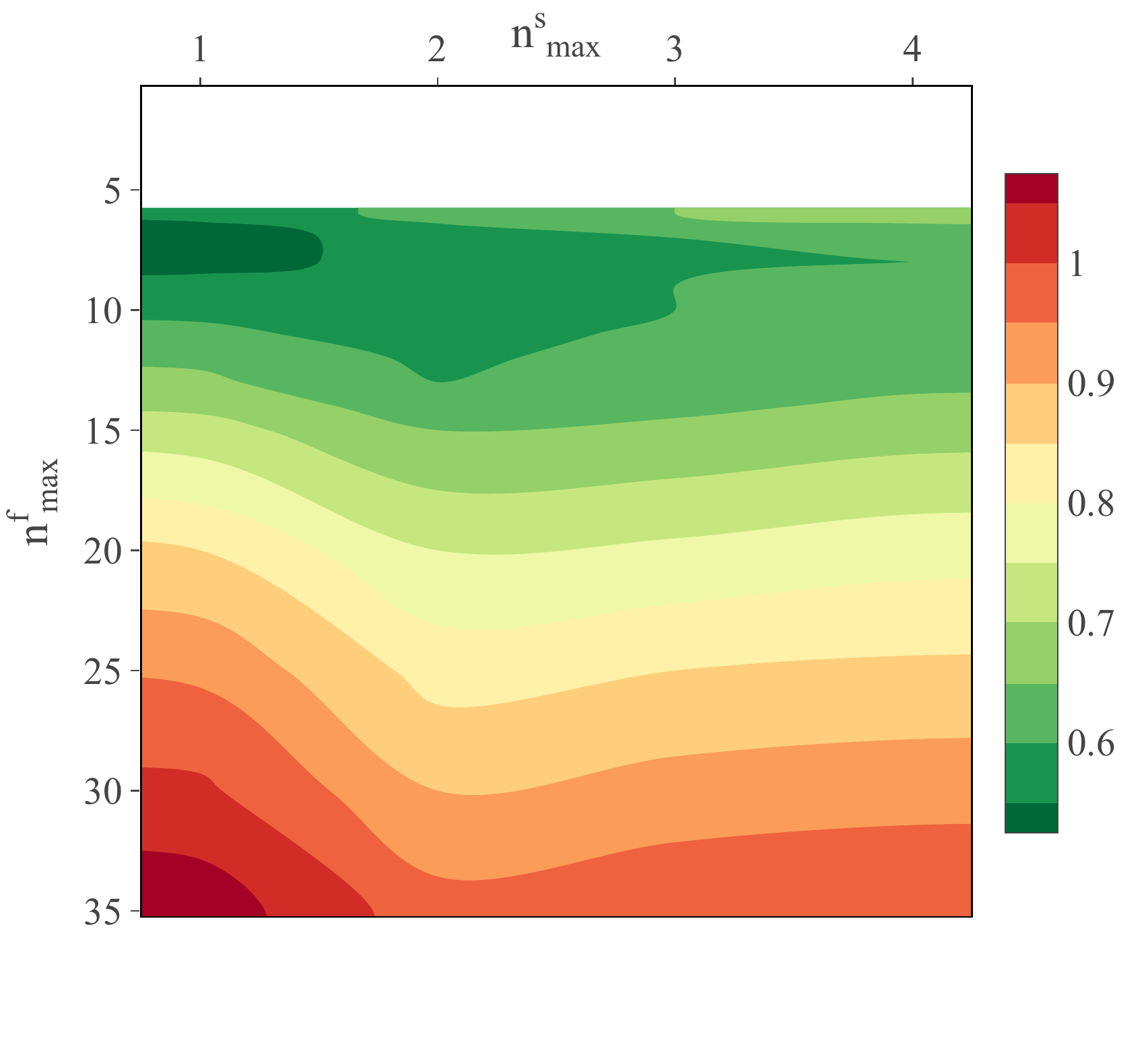} 
        \caption{$\costCoupleSum=\costIter^f=\costIter^s=1$: running only a few \subproblemIter s per call reduces the sum of iterations $\totalCoupleIter+\totalProblemIter^f+\totalProblemIter^s$, which is proportional to the overall cost for this scenario.}
    \end{subfigure}
    \begin{subfigure}[t]{.33\textwidth}
        \includegraphics[width=\textwidth,trim={0 50 10 5},clip]{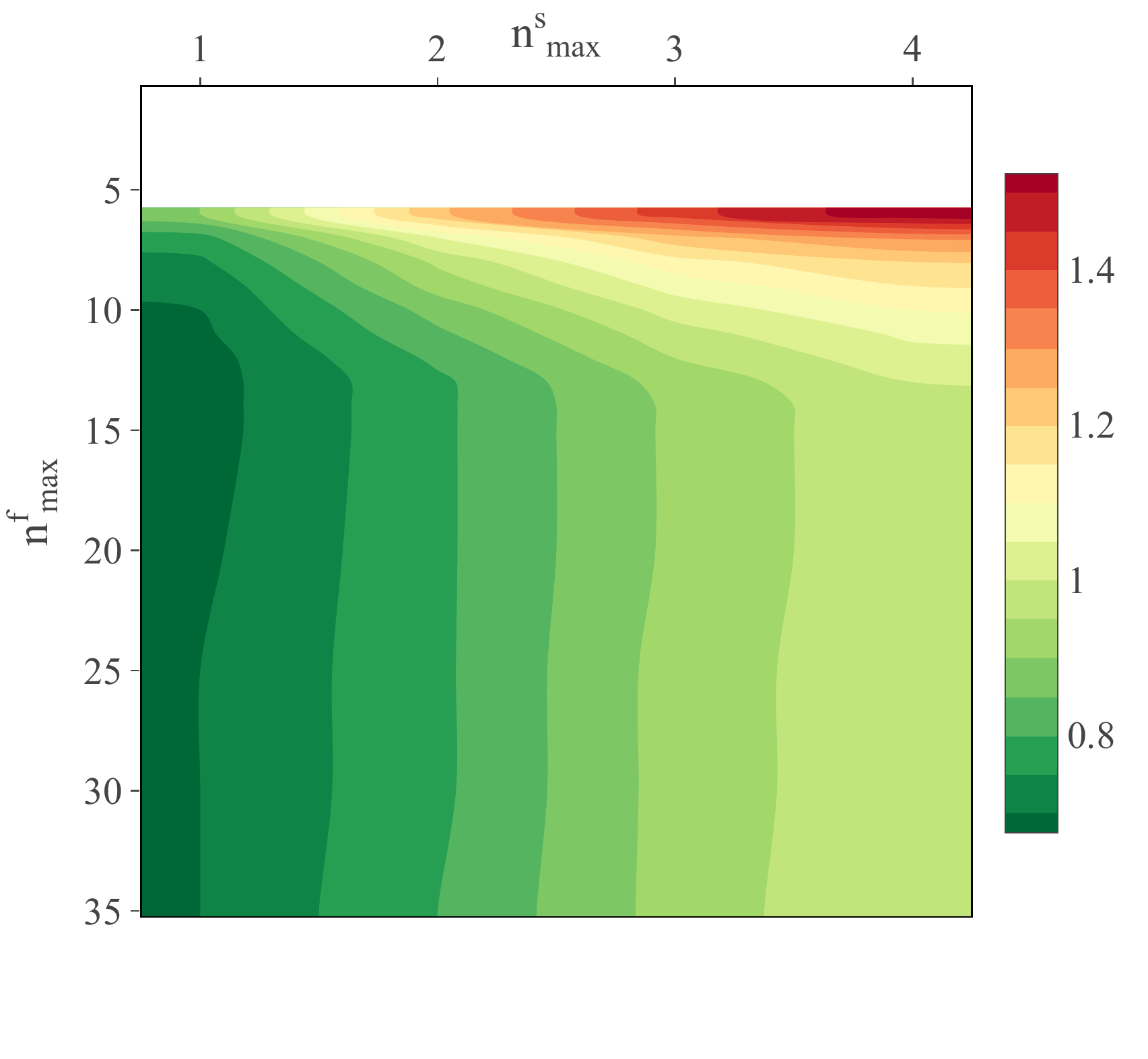}
    	\caption{$\costCoupleSum=\costIter^s=1$, $\costIter^f=0.01$: in this scenario the flow iterations are very inexpensive, so that a high $\iterPerCall^f$ and a small $\iterPerCall^s$ are most efficient.}
    \end{subfigure}
    \begin{subfigure}[t]{.33\textwidth}
        \includegraphics[width=\textwidth,trim={0 50 10 5},clip]{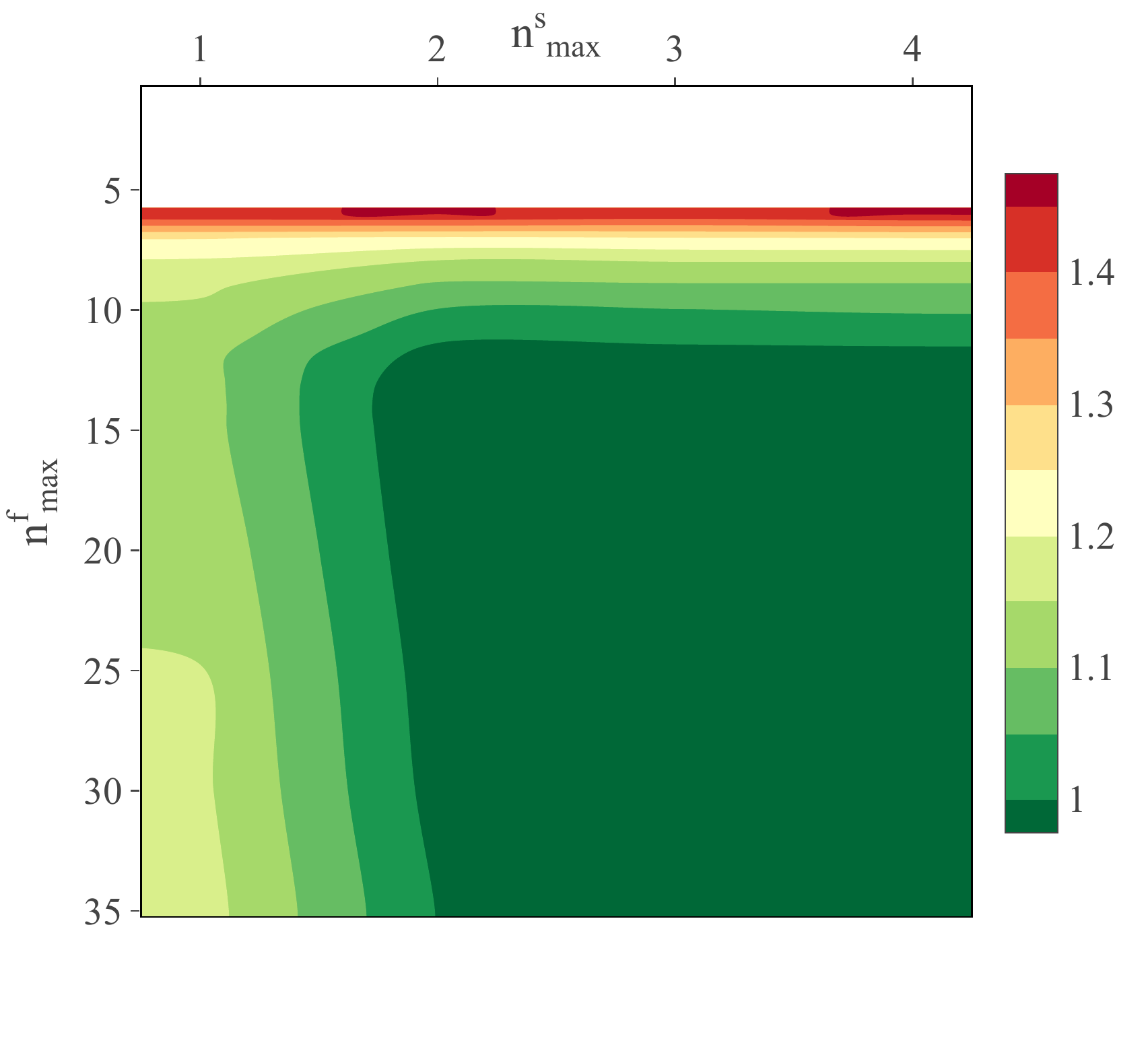}
    	\caption{$\costCoupleSum=120$, $\costIter^f=\costIter^s=1$: if the cost per coupling step becomes dominant, 
        it is best to run more \subproblemIter s per call since it reduces $\totalCoupleIter$. Thus the potential speed-up decreases.}
    \end{subfigure}
    \caption{Contour plots of different cost distributions for the cavity case. The underlying approximate data is created by using the iteration counts obtained for the FV-FE framework and artificially setting the cost per coupling and \subproblemIter.
    }
    \label{fig:ContourPlots_Varying}
\end{figure}

As an illustrative example,
\Fig{ContourPlots_Varying} uses the iteration counts obtained for the cavity case simulated in the FV-FE framework
to plot the computational cost resulting from different, artificially chosen cost factors, 
exemplarily imitating the effect of a cost-efficient flow solver (\Fig{ContourPlots_Varying}b) or a very expensive data exchange (\Fig{ContourPlots_Varying}c).


\section{Conclusion}
\label{sec:conclusion}

This work studies the influence of the nonlinear \subproblemIter s on the 
computational cost of a partitioned scheme.
Although the focus of this work is on partitioned algorithms for fluid-structure interaction problems, its findings are more widely applicable to any multi-field problem that is solved by coupling black-box solvers.
While it is common in FSI literature to consider only the total number of coupling iterations,
the \subproblemIter s are demonstrated 
to be just as crucial for the overall computational cost.
%
Consequently, 
converging in as few coupling iterations as possible 
does not necessarily minimize the run time of the partitioned simulation.
In order to be efficient, a coupling algorithm instead
should reduce the number of \subproblemIter s in the flow and solid solver too.
The cost measure proposed in this work therefore
weights these conflicting objectives with the cost of one coupling, flow or solid iteration, respectively.
%
%
Among other things, 
the cost factors depend on the considered problem, the employed solvers, and the computer hardware,
making them very case-specific.

Nevertheless,
this work sheds light on the relation between \subproblemIter s and computational cost based on typical benchmark problems and common solver types by running multiple parameter studies.
%
%
%
%
These studies systematically investigate how limiting the maximum number of \subproblemIter s per solver call
influences
the total number of coupling, flow, and solid iterations as well as the required run time.
To make sure all parameter sets yield the same solution fields,
a new convergence criterion is proposed
that relies solely on the convergence tolerances of the subproblem residuals, rather than introducing an own tolerance
for the fixed-point residual as commonly done in literature
This new way of assessing convergence can obviously also be applied in other FSI simulations.

The key finding of the parameter studies is that limiting
the \subproblemIter s per solver call
%
%
cannot only reduce the total number of \subproblemIter s,
but also significantly lower the simulation's overall computational cost. 
This raises the question of how many \subproblemIter s per solver call yield the most efficient partitioned algorithm.

Unfortunately, 
this question
is as non-trivial as it is important.
Although a definite answer is yet to be found and is expected to be case-specific,
the results discussed in this work demonstrate that iterating to full convergence in every solver call 
typically causes a computational overhead.
Instead, running only a few \subproblemIter s in each coupling step can be much more efficient. 
For the four parameter studies, the maximum 
reduction in computational cost was above $20\,\%$.
Depending on the coupling technique, even higher savings are possible, as indicated by the studies without reuse of past time step data in the IQN for which the cost was almost halved.
The ideal choice may be difficult or even impossible to determine a priori, 
but the results show that the optimum is rather flat.
For example, limiting the Newton iterations per solver call to $2$ and the fixed-point iterations to $12$
led to a speed up of $12\,\%$ to $24\,\%$ for all test cases and frameworks investigated.
While these performance gains might be no quantum leap, it should be noted that they come without any additional effort, cost, or code changes.

Furthermore, the limits imposed on the number of \subproblemIter s per solver call were 
kept fixed here for the sake of cleaner parameter studies that are easier to interpret. 
Dynamically adapting $\iterPerCall^f$ and $\iterPerCall^s$
in a ``smart'' way,
for example based on some quality measure for the coupling data, 
has the potential to further shorten the run time.\\


All in all, the primary goal of this work is to raise awareness
of a gap in current literature on FSI and other coupled problems,
%
concerning
the number of \subproblemIter s performed per solver call 
and how it influences the computational cost of a partitioned algorithm.
In addition, this work lays the groundwork
for potential future research bridging this gap.

\section*{Acknowledgement}
T. Spenke gratefully acknowledges the computing time granted by the JARA Vergabegremium and provided on the JARA Partition part of the supercomputer CLAIX at RWTH Aachen University.

N. Delaissé gratefully acknowledges the funding received from the Research Foundation - Flanders (FWO).

\appendix
\section{Additional regression measures}
\label{app:fullRegression}
The regression method to obtain the cost factors is explained in \Sec{regression}, illustrated by the simulations of the flexible tube with the FV-FE. The accuracy of the various regressions for the other parameter studies is summarized in \Tab{regressionAccuracy}.

Besides the RRMSE, the absolute variant is also given, namely the root mean square error (RMSE)
\begin{equation}
    \mathrm{RMSE} = \sqrt{\frac{\sum^\numCalc\abs{
    \timing^f-(\totalProblemIter^c \cdot \costFix + \totalProblemIter^f \cdot \costIter)
    }^2}{m}}.
    \label{equ:RMSE}
\end{equation}

\begin{table}[ph]
	\caption{Quality measures for the coupling, flow, and solid regression for all parameters studies, as well as deviation between the actual and equivalent run times determined with the new measure $\costSimulation$ and with the measure used in literature, i.e., only regarding the number of coupling iterations.}
    \label{tab:regressionAccuracy}
    \begin{center}
        \smallIfElsevier 
        \label{tab:RRMSE_Measures}
		\begin{tabular}{l|*{6}{l}|ll|ll}
			& \multicolumn{2}{c}{Coupling} & \multicolumn{2}{c}{Flow} & \multicolumn{2}{c|}{Solid} & \multicolumn{2}{c|}{New equivalent time} & \multicolumn{2}{c}{Measure in literature} \\
            & RMSE & RRMSE & RMSE & RRMSE & RMSE & RRMSE & MAPE & maxAPE & MAPE & maxAPE \\
            \hline
			\textbf{Lid-driven cavity} & & & & & & & & & & \\
			$\quad$ FE-FE     & \qty{0.2}{\second} & 7.88 \% 
                              & \qty{1.7}{\second} & 0.44 \%  
                              & \qty{0.3}{\second} & 1.26 \%  
                              & 0.33 \% & 0.89 \% 
                              & 22.77 \% & 34.32 \% \\
            $\quad$ FV-FE     & \qty{13.8}{\second} & 11.95 \% 
                              & \qty{468.7}{\second} & 5.60 \%  
                              & \qty{6.3}{\second} & 3.74 \%  
                              & 3.86 \% & 15.76 \%
                              & 12.29 \% & 36.20 \% \\
            \textbf{Flexible tube} & & & & & & & & & & \\
			$\quad$ FE-FE     & \qty{2.8}{\second} & 1.86 \% 
                              & \qty{7.3}{\second} & 0.32 \%  
                              & \qty{0.8}{\second} & 0.26 \%  
                              &  0.22 \% &  0.58 \% 
                              & 12.70 \% & 33.32 \% \\
			$\quad$ FV-FE     & \qty{2.7}{\second} & 2.65 \% 
                              & \qty{233.2}{\second} & 7.43 \%  
                              & \qty{27.7}{\second} & 3.44 \%  
                              & 3.80 \% & 14.62 \%
                              & 13.58 \% & 52.69 \% \\
		\end{tabular}
	\end{center}
\end{table}

\section{Additional results}
\label{app:fullResults}
To increase readability, the discussion, tables, and figures of \Sec{results} focused on one trend at a time.
For the sake of completeness, this appendix groups all findings in a single table for each parameter study, listing the equivalent time measure, the number of coupling iterations, and the number of \subproblemIter s for every run within a study.
\Tab{FECavity} and \Tab{FVCavity} list the results of the lid-driven cavity case for the FE-FE and FV-FE framework, respectively.
Similarly, the flexible tube case with the FE-FE setup is covered by \Tab{FETube}, and its FV-FE analogue by \Tab{FVTube}.

\renewcommand{\minval}{0.777}
\begin{table}
\centering
\caption{Lid-driven cavity case with the FE-FE framework.
The row and column header contain the maximal number of \subproblemIter s for the flow and solid solver, $n_{max}^{f}$ and $n_{max}^{s}$, respectively. For each run, the normalized \textbf{equivalent time} is given, as well as the number of \cIt{coupling iterations}, \fIt{flow solver iterations}, and \sIt{solid solver iterations}.}
\label{tab:FECavity}
\begin{tabular}{cc|cc|cc|cc|cc}
\multicolumn{2}{c}{}    & \multicolumn{8}{c}{Newton iterations per coupling iteration - Structural solver} \\
\multicolumn{2}{c}{}    & \multicolumn{2}{c|}{1} & \multicolumn{2}{c|}{2} & \multicolumn{2}{c|}{3} & \multicolumn{2}{c}{$\infty$} \\ \cline{3-10}
\multirow{8}{*}{\rotatebox[origin=c]{90}{\makecell{Newton iterations \\ per coupling iteration \\ - Flow solver}}}
&                       & \eqTime{1.07} & \cIt{8554}& \eqTime{1.15} & \cIt{8876}& \eqTime{1.12} & \cIt{8542}& \eqTime{1.12} & \cIt{8499}\\
& \multirow{-2}{*}{1}   & \fIt{8554} & \sIt{8554}& \fIt{8876} & \sIt{17043}& \fIt{8542} & \sIt{18894}& \fIt{8499} & \sIt{19662}\\ \cline{2-10}
&                       & \eqTime{0.91} & \cIt{5066}& \eqTime{0.78} & \cIt{4304}& \eqTime{0.78} & \cIt{4293}& \eqTime{0.79} & \cIt{4299}\\
& \multirow{-2}{*}{2}   & \fIt{9014} & \sIt{5066}& \fIt{7470} & \sIt{7908}& \fIt{7452} & \sIt{9294}& \fIt{7473} & \sIt{10183}\\ \cline{2-10}
&                       & \eqTime{1.05} & \cIt{4985}& \eqTime{0.83} & \cIt{3989}& \eqTime{0.84} & \cIt{4014}& \eqTime{0.84} & \cIt{4002}\\
& \multirow{-2}{*}{3}   & \fIt{11098} & \sIt{4985}& \fIt{8489} & \sIt{7278}& \fIt{8515} & \sIt{8729}& \fIt{8458} & \sIt{9579}\\ \cline{2-10}
&                       & \eqTime{1.21} & \cIt{5004}& \eqTime{0.99} & \cIt{4007}& \eqTime{0.99} & \cIt{3982}& \eqTime{1.00} & \cIt{3974}\\
& \multirow{-2}{*}{$\infty$}   & \fIt{13487} & \sIt{5004}& \fIt{10853} & \sIt{7314}& \fIt{10761} & \sIt{8667}& \fIt{10786} & \sIt{9523}
\end{tabular}
\end{table}

\renewcommand{\minval}{0.759}
\begin{table}
\centering
\caption{Lid-driven cavity case with the FV-FE framework.
The row and column header contain the maximal number of \subproblemIter s for the flow and solid solver, $n_{max}^{f}$ and $n_{max}^{s}$, respectively. For each run, the normalized \textbf{equivalent time} is given, as well as the number of \cIt{coupling iterations}, \fIt{flow solver iterations}, and \sIt{solid solver iterations}. A missing value indicates that the coupling did not converge.}
\label{tab:FVCavity}
\begin{tabular}{cc|cc|cc|cc|cc}
\multicolumn{2}{c}{}    & \multicolumn{8}{c}{Newton iterations per coupling iteration - Structural solver} \\
\multicolumn{2}{c}{}    & \multicolumn{2}{c|}{1} & \multicolumn{2}{c|}{2} & \multicolumn{2}{c|}{3} & \multicolumn{2}{c}{$\infty$} \\ \cline{3-10}
\multirow{30}{*}{\rotatebox[origin=c]{90}{Fixed-point iterations per coupling iteration - Flow solver}}
&                       & \eqTime{-} & \cIt{-}& \eqTime{-} & \cIt{-}& \eqTime{-} & \cIt{-}& \eqTime{-} & \cIt{-}\\
& \multirow{-2}{*}{5}   & \fIt{-} & \sIt{-}& \fIt{-} & \sIt{-}& \fIt{-} & \sIt{-}& \fIt{-} & \sIt{-}\\ \cline{2-10}
&                       & \eqTime{0.96} & \cIt{7209}& \eqTime{0.97} & \cIt{7204}& \eqTime{0.95} & \cIt{7111}& \eqTime{0.96} & \cIt{7176}\\
& \multirow{-2}{*}{6}   & \fIt{37477} & \sIt{7209}& \fIt{37581} & \sIt{12900}& \fIt{36929} & \sIt{16148}& \fIt{37399} & \sIt{17802}\\ \cline{2-10}
&                       & \eqTime{0.86} & \cIt{6218}& \eqTime{0.85} & \cIt{6126}& \eqTime{0.85} & \cIt{6132}& \eqTime{0.85} & \cIt{6138}\\
& \multirow{-2}{*}{7}   & \fIt{36473} & \sIt{6218}& \fIt{35965} & \sIt{10865}& \fIt{36001} & \sIt{13791}& \fIt{36068} & \sIt{15110}\\ \cline{2-10}
&                       & \eqTime{0.84} & \cIt{5882}& \eqTime{0.81} & \cIt{5613}& \eqTime{0.81} & \cIt{5618}& \eqTime{0.81} & \cIt{5608}\\
& \multirow{-2}{*}{8}   & \fIt{38449} & \sIt{5882}& \fIt{36629} & \sIt{9926}& \fIt{36706} & \sIt{12648}& \fIt{36590} & \sIt{13896}\\ \cline{2-10}
&                       & \eqTime{0.84} & \cIt{5692}& \eqTime{0.79} & \cIt{5328}& \eqTime{0.79} & \cIt{5318}& \eqTime{0.79} & \cIt{5310}\\
& \multirow{-2}{*}{9}   & \fIt{40841} & \sIt{5692}& \fIt{38382} & \sIt{9443}& \fIt{38163} & \sIt{11991}& \fIt{38127} & \sIt{13234}\\ \cline{2-10}
&                       & \eqTime{0.85} & \cIt{5581}& \eqTime{0.78} & \cIt{5091}& \eqTime{0.78} & \cIt{5075}& \eqTime{0.78} & \cIt{5103}\\
& \multirow{-2}{*}{10}   & \fIt{43394} & \sIt{5581}& \fIt{39559} & \sIt{8943}& \fIt{39430} & \sIt{11354}& \fIt{39631} & \sIt{12683}\\ \cline{2-10}
&                       & \eqTime{0.86} & \cIt{5485}& \eqTime{0.77} & \cIt{4906}& \eqTime{0.77} & \cIt{4885}& \eqTime{0.77} & \cIt{4888}\\
& \multirow{-2}{*}{11}   & \fIt{45919} & \sIt{5485}& \fIt{40814} & \sIt{8613}& \fIt{40662} & \sIt{10865}& \fIt{40715} & \sIt{12124}\\ \cline{2-10}
&                       & \eqTime{0.87} & \cIt{5398}& \eqTime{0.76} & \cIt{4706}& \eqTime{0.76} & \cIt{4718}& \eqTime{0.76} & \cIt{4720}\\
& \multirow{-2}{*}{12}   & \fIt{48329} & \sIt{5398}& \fIt{41953} & \sIt{8279}& \fIt{42020} & \sIt{10519}& \fIt{41978} & \sIt{11790}\\ \cline{2-10}
&                       & \eqTime{0.89} & \cIt{5387}& \eqTime{0.76} & \cIt{4607}& \eqTime{0.76} & \cIt{4602}& \eqTime{0.76} & \cIt{4588}\\
& \multirow{-2}{*}{13}   & \fIt{51042} & \sIt{5387}& \fIt{43364} & \sIt{8056}& \fIt{43320} & \sIt{10201}& \fIt{43319} & \sIt{11447}\\ \cline{2-10}
&                       & \eqTime{0.91} & \cIt{5373}& \eqTime{0.77} & \cIt{4566}& \eqTime{0.77} & \cIt{4545}& \eqTime{0.77} & \cIt{4546}\\
& \multirow{-2}{*}{14}   & \fIt{53842} & \sIt{5373}& \fIt{45288} & \sIt{7982}& \fIt{45168} & \sIt{10091}& \fIt{45386} & \sIt{11354}\\ \cline{2-10}
&                       & \eqTime{0.92} & \cIt{5352}& \eqTime{0.79} & \cIt{4558}& \eqTime{0.78} & \cIt{4534}& \eqTime{0.78} & \cIt{4534}\\
& \multirow{-2}{*}{15}   & \fIt{56336} & \sIt{5352}& \fIt{47478} & \sIt{7966}& \fIt{47241} & \sIt{10068}& \fIt{47344} & \sIt{11351}\\ \cline{2-10}
&                       & \eqTime{1.01} & \cIt{5349}& \eqTime{0.85} & \cIt{4519}& \eqTime{0.85} & \cIt{4487}& \eqTime{0.85} & \cIt{4500}\\
& \multirow{-2}{*}{20}   & \fIt{67992} & \sIt{5349}& \fIt{56902} & \sIt{7903}& \fIt{56645} & \sIt{9958}& \fIt{56739} & \sIt{11219}\\ \cline{2-10}
&                       & \eqTime{1.07} & \cIt{5357}& \eqTime{0.91} & \cIt{4511}& \eqTime{0.91} & \cIt{4498}& \eqTime{0.90} & \cIt{4483}\\
& \multirow{-2}{*}{25}   & \fIt{76627} & \sIt{5357}& \fIt{64779} & \sIt{7896}& \fIt{64833} & \sIt{9994}& \fIt{64530} & \sIt{11213}\\ \cline{2-10}
&                       & \eqTime{1.12} & \cIt{5345}& \eqTime{0.95} & \cIt{4499}& \eqTime{0.96} & \cIt{4503}& \eqTime{0.95} & \cIt{4491}\\
& \multirow{-2}{*}{30}   & \fIt{83527} & \sIt{5345}& \fIt{71215} & \sIt{7886}& \fIt{71336} & \sIt{9997}& \fIt{71008} & \sIt{11211}\\ \cline{2-10}
&                       & \eqTime{1.17} & \cIt{5348}& \eqTime{1.00} & \cIt{4523}& \eqTime{1.00} & \cIt{4492}& \eqTime{1.00} & \cIt{4509}\\
& \multirow{-2}{*}{$\infty$}   & \fIt{89966} & \sIt{5348}& \fIt{77618} & \sIt{7922}& \fIt{77344} & \sIt{9983}& \fIt{77272} & \sIt{11221}
\end{tabular}
\end{table}

\renewcommand{\minval}{0.794}
\begin{table}
\centering
\caption{Flexible tube case with the FE-FE framework.
The row and column header contain the maximal number of \subproblemIter s for the flow and solid solver, $n_{max}^{f}$ and $n_{max}^{s}$, respectively. For each run, the normalized \textbf{equivalent time} is given, as well as the number of \cIt{coupling iterations}, \fIt{flow solver iterations}, and \sIt{solid solver iterations}.}
\label{tab:FETube}
\begin{tabular}{cc|cc|cc|cc|cc}
\multicolumn{2}{c}{}    & \multicolumn{8}{c}{Newton iterations per coupling iteration - Structural solver} \\
\multicolumn{2}{c}{}    & \multicolumn{2}{c|}{1} & \multicolumn{2}{c|}{2} & \multicolumn{2}{c|}{3} & \multicolumn{2}{c}{$\infty$} \\ \cline{3-10}
\multirow{8}{*}{\rotatebox[origin=c]{90}{\makecell{Newton iterations \\ per coupling iteration \\ - Flow solver}}}
&                       & \eqTime{0.79} & \cIt{920}& \eqTime{0.80} & \cIt{863}& \eqTime{0.82} & \cIt{866}& \eqTime{0.82} & \cIt{870}\\
& \multirow{-2}{*}{1}   & \fIt{920} & \sIt{920}& \fIt{863} & \sIt{1626}& \fIt{866} & \sIt{1840}& \fIt{870} & \sIt{1879}\\ \cline{2-10}
&                       & \eqTime{0.91} & \cIt{815}& \eqTime{0.88} & \cIt{757}& \eqTime{0.89} & \cIt{758}& \eqTime{0.89} & \cIt{752}\\
& \multirow{-2}{*}{2}   & \fIt{1218} & \sIt{815}& \fIt{1109} & \sIt{1414}& \fIt{1112} & \sIt{1623}& \fIt{1107} & \sIt{1642}\\ \cline{2-10}
&                       & \eqTime{1.03} & \cIt{811}& \eqTime{0.99} & \cIt{763}& \eqTime{1.00} & \cIt{761}& \eqTime{1.00} & \cIt{755}\\
& \multirow{-2}{*}{3}   & \fIt{1461} & \sIt{811}& \fIt{1328} & \sIt{1426}& \fIt{1326} & \sIt{1628}& \fIt{1321} & \sIt{1647}\\ \cline{2-10}
&                       & \eqTime{1.03} & \cIt{811}& \eqTime{0.99} & \cIt{763}& \eqTime{1.00} & \cIt{761}& \eqTime{1.00} & \cIt{755}\\
& \multirow{-2}{*}{$\infty$}   & \fIt{1461} & \sIt{811}& \fIt{1328} & \sIt{1426}& \fIt{1326} & \sIt{1628}& \fIt{1321} & \sIt{1647}
\end{tabular}
\end{table}

\renewcommand{\minval}{0.772}
\begin{table}
\centering
\caption{Flexible tube case with the FV-FE framework.
The row and column header contain the maximal number of \subproblemIter s for the flow and solid solver, $n_{max}^{f}$ and $n_{max}^{s}$, respectively. For each run, the normalized \textbf{equivalent time} is given, as well as the number of \cIt{coupling iterations}, \fIt{flow solver iterations}, and \sIt{solid solver iterations}. A missing value indicates that the coupling did not converge.}
\label{tab:FVTube}
\begin{tabular}{cc|cc|cc|cc|cc}
\multicolumn{2}{c}{}    & \multicolumn{8}{c}{Newton iterations per coupling iteration - Structural solver} \\
\multicolumn{2}{c}{}    & \multicolumn{2}{c|}{1} & \multicolumn{2}{c|}{2} & \multicolumn{2}{c|}{3} & \multicolumn{2}{c}{$\infty$} \\ \cline{3-10}
\multirow{30}{*}{\rotatebox[origin=c]{90}{Fixed-point iterations per coupling iteration - Flow solver}}
&                       & \eqTime{1.14} & \cIt{2135}& \eqTime{-} & \cIt{-}& \eqTime{-} & \cIt{-}& \eqTime{-} & \cIt{-}\\
& \multirow{-2}{*}{7}   & \fIt{13567} & \sIt{2135}& \fIt{-} & \sIt{-}& \fIt{-} & \sIt{-}& \fIt{-} & \sIt{-}\\ \cline{2-10}
&                       & \eqTime{0.95} & \cIt{1716}& \eqTime{1.10} & \cIt{1814}& \eqTime{-} & \cIt{-}& \eqTime{1.16} & \cIt{1749}\\
& \multirow{-2}{*}{8}   & \fIt{12120} & \sIt{1716}& \fIt{13006} & \sIt{3471}& \fIt{-} & \sIt{-}& \fIt{12438} & \sIt{5006}\\ \cline{2-10}
&                       & \eqTime{0.88} & \cIt{1544}& \eqTime{0.94} & \cIt{1517}& \eqTime{0.96} & \cIt{1481}& \eqTime{1.06} & \cIt{1586}\\
& \multirow{-2}{*}{9}   & \fIt{11896} & \sIt{1544}& \fIt{11838} & \sIt{2847}& \fIt{11581} & \sIt{3492}& \fIt{12381} & \sIt{4358}\\ \cline{2-10}
&                       & \eqTime{0.82} & \cIt{1400}& \eqTime{0.86} & \cIt{1348}& \eqTime{0.92} & \cIt{1387}& \eqTime{0.94} & \cIt{1385}\\
& \multirow{-2}{*}{10}   & \fIt{11758} & \sIt{1400}& \fIt{11351} & \sIt{2488}& \fIt{11722} & \sIt{3225}& \fIt{11667} & \sIt{3646}\\ \cline{2-10}
&                       & \eqTime{0.81} & \cIt{1355}& \eqTime{0.81} & \cIt{1245}& \eqTime{0.85} & \cIt{1254}& \eqTime{0.88} & \cIt{1279}\\
& \multirow{-2}{*}{11}   & \fIt{12242} & \sIt{1355}& \fIt{11281} & \sIt{2291}& \fIt{11430} & \sIt{2865}& \fIt{11542} & \sIt{3275}\\ \cline{2-10}
&                       & \eqTime{0.80} & \cIt{1299}& \eqTime{0.79} & \cIt{1189}& \eqTime{0.83} & \cIt{1201}& \eqTime{0.84} & \cIt{1196}\\
& \multirow{-2}{*}{12}   & \fIt{12701} & \sIt{1299}& \fIt{11457} & \sIt{2175}& \fIt{11614} & \sIt{2695}& \fIt{11584} & \sIt{3055}\\ \cline{2-10}
&                       & \eqTime{0.82} & \cIt{1299}& \eqTime{0.78} & \cIt{1154}& \eqTime{0.80} & \cIt{1143}& \eqTime{0.81} & \cIt{1133}\\
& \multirow{-2}{*}{13}   & \fIt{13455} & \sIt{1299}& \fIt{11767} & \sIt{2093}& \fIt{11717} & \sIt{2552}& \fIt{11621} & \sIt{2835}\\ \cline{2-10}
&                       & \eqTime{0.84} & \cIt{1289}& \eqTime{0.77} & \cIt{1114}& \eqTime{0.80} & \cIt{1114}& \eqTime{0.83} & \cIt{1137}\\
& \multirow{-2}{*}{14}   & \fIt{14155} & \sIt{1289}& \fIt{12101} & \sIt{2036}& \fIt{12085} & \sIt{2462}& \fIt{12278} & \sIt{2852}\\ \cline{2-10}
&                       & \eqTime{0.85} & \cIt{1278}& \eqTime{0.78} & \cIt{1104}& \eqTime{0.80} & \cIt{1101}& \eqTime{0.80} & \cIt{1076}\\
& \multirow{-2}{*}{15}   & \fIt{14872} & \sIt{1278}& \fIt{12604} & \sIt{2013}& \fIt{12540} & \sIt{2430}& \fIt{12156} & \sIt{2658}\\ \cline{2-10}
&                       & \eqTime{0.93} & \cIt{1264}& \eqTime{0.82} & \cIt{1062}& \eqTime{0.83} & \cIt{1049}& \eqTime{0.84} & \cIt{1046}\\
& \multirow{-2}{*}{20}   & \fIt{18399} & \sIt{1264}& \fIt{14898} & \sIt{1914}& \fIt{14672} & \sIt{2280}& \fIt{14534} & \sIt{2570}\\ \cline{2-10}
&                       & \eqTime{1.05} & \cIt{1255}& \eqTime{0.91} & \cIt{1050}& \eqTime{0.91} & \cIt{1031}& \eqTime{0.92} & \cIt{1033}\\
& \multirow{-2}{*}{30}   & \fIt{23462} & \sIt{1255}& \fIt{18756} & \sIt{1873}& \fIt{18217} & \sIt{2233}& \fIt{18154} & \sIt{2518}\\ \cline{2-10}
&                       & \eqTime{1.12} & \cIt{1264}& \eqTime{0.98} & \cIt{1069}& \eqTime{0.98} & \cIt{1043}& \eqTime{0.97} & \cIt{1010}\\
& \multirow{-2}{*}{40}   & \fIt{26140} & \sIt{1264}& \fIt{21584} & \sIt{1922}& \fIt{21056} & \sIt{2265}& \fIt{20296} & \sIt{2486}\\ \cline{2-10}
&                       & \eqTime{1.14} & \cIt{1273}& \eqTime{1.00} & \cIt{1064}& \eqTime{0.98} & \cIt{1021}& \eqTime{1.01} & \cIt{1037}\\
& \multirow{-2}{*}{50}   & \fIt{26768} & \sIt{1273}& \fIt{22083} & \sIt{1932}& \fIt{21210} & \sIt{2206}& \fIt{21368} & \sIt{2564}\\ \cline{2-10}
&                       & \eqTime{1.14} & \cIt{1263}& \eqTime{0.99} & \cIt{1052}& \eqTime{0.98} & \cIt{1030}& \eqTime{1.00} & \cIt{1028}\\
& \multirow{-2}{*}{$\infty$}   & \fIt{27089} & \sIt{1263}& \fIt{21898} & \sIt{1892}& \fIt{21292} & \sIt{2206}& \fIt{21282} & \sIt{2543}
\end{tabular}
\end{table}

\section{Additional results without reuse}
\label{app:noReuse}

\Tab{FETubeNoReuse} and \Tab{FVTubeNoReuse} present two additional parameter studies for the FE-FE and FV-FE framework, respectively,
simulating the flexible tube case with identical settings, but without reuse of data from past time steps in the IQN Jacobian approximation ($q$=0).

\renewcommand{\minval}{0.526}
\begin{table}
\centering
\caption{Flexible tube case simulated with the FE-FE framework, but without reusing past data in the IQN Jacobian approximation. 
The formatting is the same as in the previous tables.
The cost factors determined by the regression model for this parameter study are, in seconds,
$\costIter^f=0.2177$, $\costIter^s=0.2253$ and $\costCoupleSum=0.1620$.}
\label{tab:FETubeNoReuse}
\begin{tabular}{cc|cc|cc|cc|cc}
\multicolumn{2}{c}{}    & \multicolumn{8}{c}{Newton iterations per coupling iteration - Structural solver} \\
\multicolumn{2}{c}{}    & \multicolumn{2}{c|}{1} & \multicolumn{2}{c|}{2} & \multicolumn{2}{c|}{3} & \multicolumn{2}{c}{$\infty$} \\ \cline{3-10}
\multirow{8}{*}{\rotatebox[origin=c]{90}{\makecell{Newton iterations \\ per coupling iteration \\ - Flow solver}}}
&                       & \eqTime{0.53} & \cIt{1465}& \eqTime{0.69} & \cIt{1417}& \eqTime{0.76} & \cIt{1417}& \eqTime{0.77} & \cIt{1417}\\
& \multirow{-2}{*}{1}   & \fIt{1465} & \sIt{1465}& \fIt{1417} & \sIt{2734}& \fIt{1417} & \sIt{3330}& \fIt{1417} & \sIt{3362}\\ \cline{2-10}
&                       & \eqTime{0.67} & \cIt{1463}& \eqTime{0.82} & \cIt{1416}& \eqTime{0.90} & \cIt{1416}& \eqTime{0.91} & \cIt{1416}\\
& \multirow{-2}{*}{2}   & \fIt{2573} & \sIt{1463}& \fIt{2479} & \sIt{2732}& \fIt{2479} & \sIt{3328}& \fIt{2479} & \sIt{3360}\\ \cline{2-10}
&                       & \eqTime{0.77} & \cIt{1463}& \eqTime{0.92} & \cIt{1416}& \eqTime{1.00} & \cIt{1416}& \eqTime{1.00} & \cIt{1416}\\
& \multirow{-2}{*}{3}   & \fIt{3318} & \sIt{1463}& \fIt{3206} & \sIt{2732}& \fIt{3206} & \sIt{3328}& \fIt{3206} & \sIt{3360}\\ \cline{2-10}
&                       & \eqTime{0.77} & \cIt{1463}& \eqTime{0.92} & \cIt{1416}& \eqTime{1.00} & \cIt{1416}& \eqTime{1.00} & \cIt{1416}\\
& \multirow{-2}{*}{$\infty$}   & \fIt{3318} & \sIt{1463}& \fIt{3206} & \sIt{2732}& \fIt{3206} & \sIt{3328}& \fIt{3206} & \sIt{3360}
\end{tabular}
\end{table}

\renewcommand{\minval}{0.497}
\begin{table}
\centering
\caption{Flexible tube case simulated with the FV-FE framework, but without reusing past data in the IQN Jacobian approximation. 
The formatting is the same as in the previous tables.
The cost factors determined by the regression model for this parameter study are, in seconds,
$\costIter^f=0.1139$, $\costIter^s=0.2410$, and $\costCoupleSum=1.2825$.}
\label{tab:FVTubeNoReuse}
\begin{tabular}{cc|cc|cc|cc|cc}
\multicolumn{2}{c}{}    & \multicolumn{8}{c}{Newton iterations per coupling iteration - Structural solver} \\
\multicolumn{2}{c}{}    & \multicolumn{2}{c|}{1} & \multicolumn{2}{c|}{2} & \multicolumn{2}{c|}{3} & \multicolumn{2}{c}{$\infty$} \\ \cline{3-10}
\multirow{30}{*}{\rotatebox[origin=c]{90}{Fixed-point iterations per coupling iteration - Flow solver}}
&                       & \eqTime{0.54} & \cIt{2077}& \eqTime{0.61} & \cIt{2120}& \eqTime{0.61} & \cIt{2044}& \eqTime{0.64} & \cIt{2060}\\
& \multirow{-2}{*}{5}   & \fIt{9641} & \sIt{2077}& \fIt{9815} & \sIt{3970}& \fIt{9444} & \sIt{4721}& \fIt{9511} & \sIt{5286}\\ \cline{2-10}
&                       & \eqTime{0.50} & \cIt{1834}& \eqTime{0.55} & \cIt{1827}& \eqTime{0.57} & \cIt{1827}& \eqTime{0.59} & \cIt{1828}\\
& \multirow{-2}{*}{6}   & \fIt{10097} & \sIt{1834}& \fIt{10052} & \sIt{3427}& \fIt{10050} & \sIt{4299}& \fIt{10045} & \sIt{4892}\\ \cline{2-10}
&                       & \eqTime{0.50} & \cIt{1746}& \eqTime{0.54} & \cIt{1757}& \eqTime{0.57} & \cIt{1763}& \eqTime{0.59} & \cIt{1756}\\
& \multirow{-2}{*}{7}   & \fIt{10960} & \sIt{1746}& \fIt{10973} & \sIt{3261}& \fIt{10986} & \sIt{4101}& \fIt{10965} & \sIt{4715}\\ \cline{2-10}
&                       & \eqTime{0.50} & \cIt{1691}& \eqTime{0.54} & \cIt{1681}& \eqTime{0.57} & \cIt{1685}& \eqTime{0.59} & \cIt{1685}\\
& \multirow{-2}{*}{8}   & \fIt{11881} & \sIt{1691}& \fIt{11844} & \sIt{3101}& \fIt{11880} & \sIt{3928}& \fIt{11899} & \sIt{4582}\\ \cline{2-10}
&                       & \eqTime{0.50} & \cIt{1638}& \eqTime{0.54} & \cIt{1610}& \eqTime{0.57} & \cIt{1620}& \eqTime{0.59} & \cIt{1633}\\
& \multirow{-2}{*}{9}   & \fIt{12889} & \sIt{1638}& \fIt{12767} & \sIt{2994}& \fIt{12811} & \sIt{3807}& \fIt{12841} & \sIt{4490}\\ \cline{2-10}
&                       & \eqTime{0.51} & \cIt{1595}& \eqTime{0.55} & \cIt{1578}& \eqTime{0.57} & \cIt{1583}& \eqTime{0.59} & \cIt{1593}\\
& \multirow{-2}{*}{10}   & \fIt{13941} & \sIt{1595}& \fIt{13729} & \sIt{2923}& \fIt{13758} & \sIt{3727}& \fIt{13772} & \sIt{4418}\\ \cline{2-10}
&                       & \eqTime{0.52} & \cIt{1579}& \eqTime{0.55} & \cIt{1555}& \eqTime{0.58} & \cIt{1559}& \eqTime{0.60} & \cIt{1559}\\
& \multirow{-2}{*}{11}   & \fIt{15046} & \sIt{1579}& \fIt{14713} & \sIt{2865}& \fIt{14741} & \sIt{3672}& \fIt{14741} & \sIt{4360}\\ \cline{2-10}
&                       & \eqTime{0.54} & \cIt{1572}& \eqTime{0.56} & \cIt{1531}& \eqTime{0.59} & \cIt{1523}& \eqTime{0.61} & \cIt{1520}\\
& \multirow{-2}{*}{12}   & \fIt{16166} & \sIt{1572}& \fIt{15756} & \sIt{2820}& \fIt{15747} & \sIt{3608}& \fIt{15745} & \sIt{4302}\\ \cline{2-10}
&                       & \eqTime{0.55} & \cIt{1564}& \eqTime{0.58} & \cIt{1520}& \eqTime{0.60} & \cIt{1515}& \eqTime{0.62} & \cIt{1515}\\
& \multirow{-2}{*}{13}   & \fIt{17309} & \sIt{1564}& \fIt{16836} & \sIt{2787}& \fIt{16814} & \sIt{3577}& \fIt{16815} & \sIt{4280}\\ \cline{2-10}
&                       & \eqTime{0.57} & \cIt{1562}& \eqTime{0.59} & \cIt{1511}& \eqTime{0.61} & \cIt{1509}& \eqTime{0.63} & \cIt{1508}\\
& \multirow{-2}{*}{14}   & \fIt{18463} & \sIt{1562}& \fIt{17905} & \sIt{2767}& \fIt{17872} & \sIt{3565}& \fIt{17863} & \sIt{4267}\\ \cline{2-10}
&                       & \eqTime{0.59} & \cIt{1563}& \eqTime{0.60} & \cIt{1502}& \eqTime{0.63} & \cIt{1500}& \eqTime{0.65} & \cIt{1500}\\
& \multirow{-2}{*}{15}   & \fIt{19615} & \sIt{1563}& \fIt{18969} & \sIt{2746}& \fIt{18915} & \sIt{3541}& \fIt{18914} & \sIt{4246}\\ \cline{2-10}
&                       & \eqTime{0.67} & \cIt{1564}& \eqTime{0.67} & \cIt{1483}& \eqTime{0.70} & \cIt{1474}& \eqTime{0.72} & \cIt{1470}\\
& \multirow{-2}{*}{20}   & \fIt{25033} & \sIt{1564}& \fIt{24210} & \sIt{2684}& \fIt{24124} & \sIt{3471}& \fIt{24107} & \sIt{4178}\\ \cline{2-10}
&                       & \eqTime{0.80} & \cIt{1564}& \eqTime{0.81} & \cIt{1479}& \eqTime{0.83} & \cIt{1472}& \eqTime{0.85} & \cIt{1472}\\
& \multirow{-2}{*}{30}   & \fIt{34388} & \sIt{1564}& \fIt{33540} & \sIt{2676}& \fIt{33450} & \sIt{3461}& \fIt{33440} & \sIt{4170}\\ \cline{2-10}
&                       & \eqTime{0.90} & \cIt{1564}& \eqTime{0.91} & \cIt{1479}& \eqTime{0.93} & \cIt{1472}& \eqTime{0.95} & \cIt{1472}\\
& \multirow{-2}{*}{40}   & \fIt{41504} & \sIt{1564}& \fIt{40657} & \sIt{2675}& \fIt{40565} & \sIt{3461}& \fIt{40556} & \sIt{4170}\\ \cline{2-10}
&                       & \eqTime{0.93} & \cIt{1564}& \eqTime{0.94} & \cIt{1479}& \eqTime{0.96} & \cIt{1472}& \eqTime{0.99} & \cIt{1472}\\
& \multirow{-2}{*}{50}   & \fIt{43672} & \sIt{1564}& \fIt{42838} & \sIt{2675}& \fIt{42749} & \sIt{3461}& \fIt{42741} & \sIt{4170}\\ \cline{2-10}
&                       & \eqTime{0.93} & \cIt{1564}& \eqTime{0.96} & \cIt{1479}& \eqTime{0.98} & \cIt{1472}& \eqTime{1.00} & \cIt{1472}\\
& \multirow{-2}{*}{$\infty$}   & \fIt{43677} & \sIt{1564}& \fIt{43793} & \sIt{2675}& \fIt{43704} & \sIt{3461}& \fIt{43697} & \sIt{4170}
\end{tabular}
\end{table}


\bibliography{journals.bib, references.bib}

\end{document}